\documentclass[aip,pof,preprint]{revtex4-1}
\usepackage{amsmath}
\usepackage{amssymb}
\usepackage{graphicx}
\usepackage{enumerate}
\usepackage{subfig}
\usepackage{comment}
\usepackage[]{natbib}
\usepackage[colorlinks=true,citecolor=blue]{hyperref}

\include{./revtex4-1}

\newcommand{\id}{\mbox{d}}

\newcommand{\vc}{\mathbf}

\makeatother

\begin{document}

\title{{How coherent are the vortices of two-dimensional
turbulence?}}
\author{{Mohammad Farazmand}}
\affiliation{{Department of Mathematics, ETH Zurich, Switzerland}}
\affiliation{{Institute of Mechanical Systems, ETH Zurich, Switzerland}}
\author{{George Haller}}
\affiliation{{Institute of Mechanical Systems, ETH Zurich, Switzerland}}
\date{{\today}}

\begin{abstract}
{We use recent developments in the theory of finite-time
dynamical systems to locate the material boundaries of coherent vortices
objectively in two-dimensional Navier--Stokes turbulence. We show
that these boundaries are optimal in the sense that any closed curve
in their exterior will lose coherence under material advection. Through
a detailed comparison, we find that other available Eulerian and Lagrangian
techniques significantly underestimate the size of each coherent vortex. }
\end{abstract}
\maketitle

\section{Introduction}
{Coherent vortices are persistent patches of rotating
fluid that are observed in experimentally and numerically generated
two-dimensional turbulence \citep{provenzale99,kellay2002,provenzale2008coherent,2dturb_rev}.
As opposed to a typical closed material line, the boundary of a coherent
vortex is envisioned to preserve its overall shape without substantial
stretching, folding or filamentation \cite{fiedler1988,bhEddy,bolt_shapeCoh}.
While intuitive and simple, this material view on vortices has proven
surprisingly challenging to implement in detecting vortex boundaries
\citep{frischbook,GFDfund}.}

The formation and evolution of coherent fluid blobs is part of the
material response of the fluid to external effects. By a classic axiom
of continuum mechanics \cite{gurtin1981}, this material response
should be objective, i.e., invariant with respect to time-dependent rotations and
translations of the frame of the observer. Yet vorticity, the main
diagnostic for structure identification in fluid mechanics, is not
objective: it changes in coordinate systems rotating relative to each
other, thus giving conflicting vortex definitions in different frames.
Consequently, there is no well-justified threshold over which vorticity
should necessarily mark a vortex.

{To address this issue, a number of alternative Eulerian
diagnostics have been proposed for vortex detection (see \citet{vortexIdent}
and \citet{objVortex}, for a review). For instance, the Okubo--Weiss
(OW) criterion \citep{okubo,weiss_okubo} identifies vortices as regions
where vorticity dominates strain. The $Q$-criterion offers a three-dimensional
version of this principle}\citep{chong1990}.
In later work, \citet{hua1998} also include accelerations in the
strain-vorticity comparison. Unfortunately, all these instantaneous
diagnostics still lack objectivity, as well as a clearly derived mathematical
connection to sustained material coherence. As a consequence, vortex
boundaries suggested by instantaneous Eulerian diagnostics tend to
lose their coherence rapidly under advection in unsteady flows \cite{javier_agulhas}.

{A recent development in the theory of finite-time
dynamical systems \citep{bhEddy} offers an objective (frame-independent)
and threshold-free Lagrangian approach to the identification of coherent
vortices in two-dimensional flows. Specifically, \citet{bhEddy} show
that coherent (non-filamenting) material lines are necessarily stationary
curves of an appropriately defined Lagrangian strain-energy functional.
They solve this variational problem explicitly to uncover vortex boundaries
as outermost limit cycles of a vector field derived from the invariants
of the Cauchy--Green strain tensor. \citet{bhEddy} demonstrate the
efficacy of this approach by extracting maximally coherent Agulhas
rings from satellite-derived oceanic surface velocities.}

{Here, we use this method to detect the optimal boundaries
of coherent vortices in a direct numerical simulation of Navier--Stokes
turbulence. We also carry out a detailed comparison with alternative
Eulerian and Lagrangian techniques. This comparison reveals that the
coherent vortices that survive for long times are significantly larger
than previously thought.}

\section{Preliminaries}\label{sec:theory}
\subsection{Set-up}
{Let $u(x,t)$ be a two-dimensional velocity field,
defined over positions $x$ in an open domain $U\subset\mathbb{R}^{2}$
and times $t$ ranging though a finite interval $I=[a,b]$. We assume
that $u(x,t)$ is a continuously differentiable function of its arguments.
The motion of passive fluid particles under such a velocity field
is governed by the differential equation 
\begin{equation}
\dot{x}=u(x,t),\label{eq:dynsys}
\end{equation}
where $x(t;t_{0},x_{0})$ is the position of a particle at time $t$
whose initial position at time $t_{0}$ is $x_{0}\in U$. For the
fixed time interval $I$, the dynamical system (\ref{eq:dynsys})
defines the flow map 
\begin{align}
F:\  & U\rightarrow U,\nonumber \\
 & x_{a}\mapsto x_{b},
\end{align}
which takes an initial condition $x_{a}$ to its time-$b$ position
$x_{b}=F(x_{a}):=x(b;a,x_{a})$. We recall form the classic theory
of ordinary differential equations that the flow map $F$ is as smooth
as the underlying velocity field $u$. \citep{arnold78}}

\subsection{Coherence principle}
{A typical set of fluid particles deforms significantly
as advected under the flow map $F$, provided that the advection time
$b-a$ is at least of the order of a few eddy turn-over times in a
turbulent flow \citep{aref}. One may seek coherent material vortices
as atypical sets of fluid trajectories that defy this trend by preserving
their overall shape. These shapes are necessarily bounded by closed
material lines that rotate and translate, but otherwise show no appreciable
stretching or folding.}

\begin{figure}
\centering\includegraphics[width=0.8\textwidth]{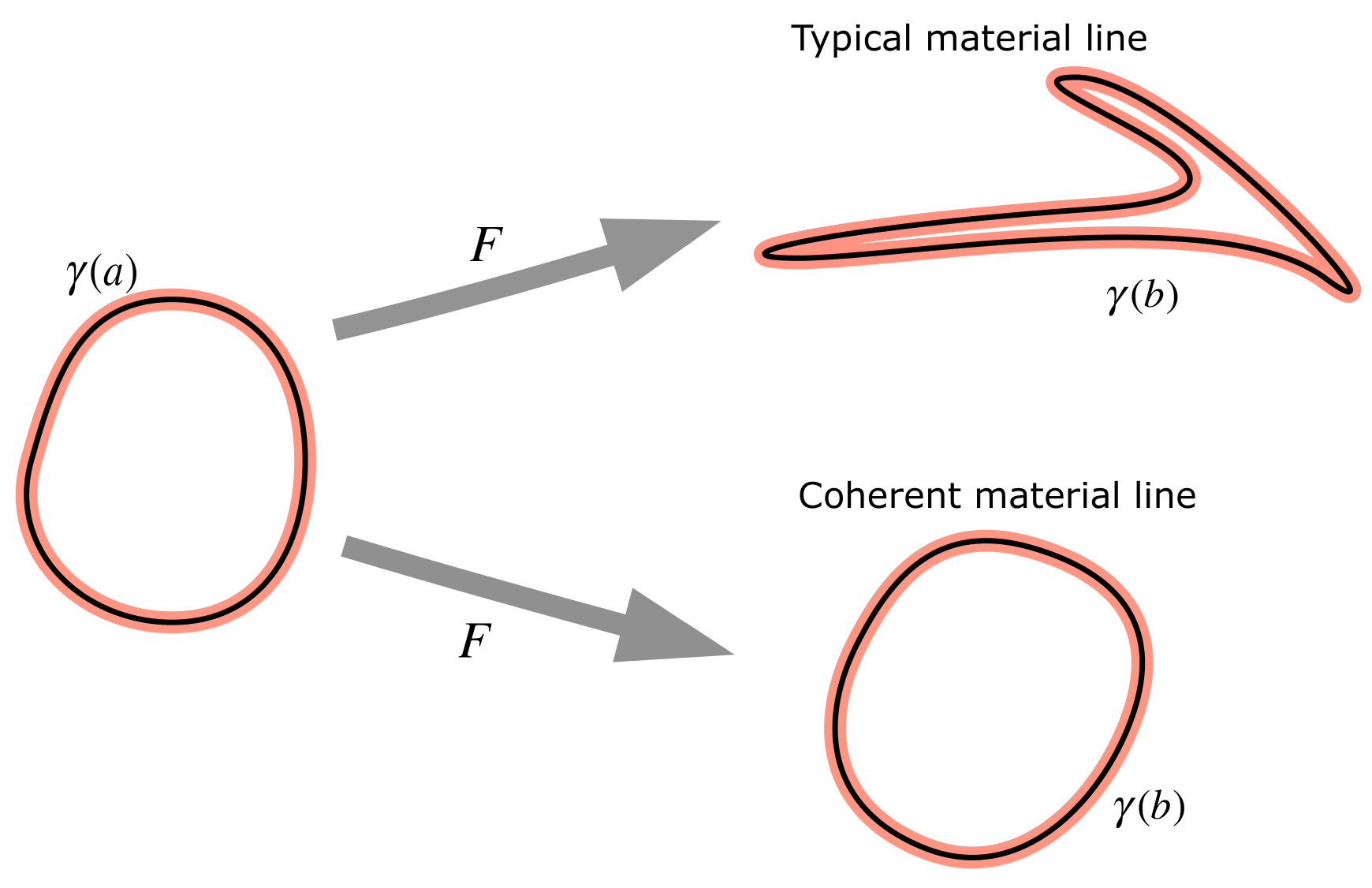}
\protect\caption{Deformation of a typical material line $\gamma$ versus a coherent
material line as advected under the flow map $F$ from time $t=a$
to $t=b$. No leading-order variation in the averaged stretching is
observed in the material belt (red) around a coherent material line.}
\label{fig:bhBoundary_sketch}
\end{figure}

{\citet{bhEddy} seek Lagrangian vortex boundaries
as closed material lines across which the averaged material straining
shows no leading-order variability. Specifically, a thin material
belt around a typical material line $\gamma$ experiences visible
inhomogeneity in straining under advection. A thin material belt around
a coherent material line, however, does not exhibit leading-order
inhomogeneity in its straining (see figure \ref{fig:bhBoundary_sketch}).}

To formulate this observation mathematically, we let $\gamma$
be a closed material line over the time interval $[a,b]$, and let
$r:s\mapsto r(s)$, with $s\in[0,\sigma]$, be a parametrization of
$\gamma$ at the initial time $t=a$. The averaged tangential strain
along $\gamma$, computed between the times $a$ and $b$ is then given by 
\begin{equation}
Q(\gamma)=\frac{1}{\sigma}\int_{0}^{\sigma}\frac{\sqrt{\langle r'(s),C(r(s))r'(s)\rangle}}{\sqrt{\langle r'(s),r'(s)\rangle}}\, ds,\label{eq:avLagStr}
\end{equation}
where the right Cauchy--Green strain tensor\citep{truesdell}
$C=DF^{\top}DF$ is defined in terms of the Jacobian $DF$ of the
flow map; the symbol $\top$ denotes matrix transposition;
prime denotes differentiation with respect to
the parameter $s$; $\langle\cdot,\cdot\rangle$ denotes the Euclidean
inner product. The integrand in equation \eqref{eq:avLagStr} represents
the pointwise tangential strain (see \citet{bhEddy} for details).

Consider a small perturbation to $\gamma$ given
by $\gamma+\epsilon h$, where $0<\epsilon\ll1$ and $h:[0,\sigma]\rightarrow\mathbb{R}^{2}$
is a $\sigma$-periodic, $\mathcal{O}(1)$ vector field orthogonal
to $\gamma$. The perturbation $\gamma+\epsilon h$ represents the
thin material belt of figure \ref{fig:bhBoundary_sketch}. For a typical
material line, we have $Q(\gamma+\epsilon h)=Q(\gamma)+\mathcal{O}(\epsilon)$
owing to the smoothness of the flow map $F$. That is, $\mathcal{O}(\epsilon)$-perturbations
to the material line $\gamma$ lead to a $\mathcal{O}(\epsilon)$-perturbation
in the averaged tangential strain $Q$. \citet{bhEddy} argue that
for a thin material belt centered on $\gamma$ to remain coherent,
the belt should not show a leading-order change with respect to $\epsilon$
in its averaged straining. In other words, $Q(\gamma+\epsilon h)=Q(\gamma)+\mathcal{O}(\epsilon^{2})$
should hold for $\gamma$, or equivalently, the first variation of
$Q$ should vanish on $\gamma$, i.e., $\delta Q(\gamma)=0$.

{The Euler-Lagrange equations arising from the condition
$\delta Q(\gamma)=0$ are too complicated to yield insight. \citet{bhEddy}
show, however, that a material line satisfies $\delta Q(\gamma)=0$
if and only if it satisfies the pointwise condition 
\begin{equation}
\langle r'(s),E_{\lambda}(r(s))r'(s)\rangle=0,\label{eq:null}
\end{equation}
for some constant $\lambda>0$, with the generalized }{\emph{Green--Lagrange
strain tensor}}{{} $E_{\lambda}$ defined as 
\begin{equation}
E_{\lambda}=\frac{1}{2}[C-\lambda^{2}I],\label{eq:strain_energy}
\end{equation}
where $I$ is the two-by-two identity matrix.}

{The implicit family of differential equations \eqref{eq:null}
is equivalent to two families of explicit differential equations of
the form 
\begin{equation}
r'=\eta_{\lambda}^{\pm}(r):=\sqrt{\frac{\lambda_{2}(r)-\lambda^{2}}{\lambda_{2}(r)-\lambda_{1}(r)}}\xi_{1}(r)\pm\sqrt{\frac{\lambda^{2}-\lambda_{1}(r)}{\lambda_{2}(r)-\lambda_{1}(r)}}\xi_{2}(r),\label{eq:eta}
\end{equation}
where $0<\lambda_{1}\leq\lambda_{2}$ are eigenvalues of $C$ and
$\{\xi_{1},\xi_{2}\}$ are their corresponding orthogonal eigenvectors
\cite{bhEddy}. In an incompressible flow, we have $\lambda_{1}\lambda_{2}=1$
(see, e.g., \citet{arnold78}).}

The vectors $\eta_{\lambda}^{+}$ and $\eta_{\lambda}^{-}$ are one-parameter
families of vector fields with $\lambda$ acting as the parameter.
In an incompressible flow, we have $\lambda_{2}\geq1$ and $\lambda_{1}\leq1$.
Therefore, for $\lambda=1$, $\eta_{\lambda}^{\pm}$ are well-defined
real vector fields over the entire physical domain $U$. For $\lambda\neq1$,
the vector fields $\eta_{\lambda}^{\pm}$ are only defined over a
subset $U_{\lambda}\subset U$ where $\lambda_{2}\geq\lambda^{2}$
and $\lambda_{1}\leq\lambda^{2}$. The trajectories of $\eta_{\lambda}^{\pm}$
can be computed over $U_{\lambda}$. We refer to these trajectories
as $\lambda$-stretching material lines (or \emph{$\lambda$-lines}, for short)

\subsection{Lagrangian vortex boundaries and $\lambda$-lines}
Here, we recall from \citet{bhEddy} some properties
of the $\lambda$-lines (i.e., trajectories of \eqref{eq:eta}), that
are relevant for Lagrangian coherent vortex detection:

(i) \emph{Uniform stretching}:
$\lambda$-lines stretch uniformly by a factor of $\lambda$ as advected
under the flow map $F$. To quantify this statement, let $\gamma_{a}$
be the time-$a$ position of a $\lambda$-line parametrized by $r:s\mapsto r(s)$.
Since $\gamma_{a}$ is a $\lambda$-line, we have $r'(s)\parallel\eta_{\lambda}^{\pm}(r(s))$.
Its time-$b$ position $\gamma_{b}$ will be parametrized by $F\circ r:s\mapsto F(r(s))$,
whose tangent vector is given by $DF(r(s))r'(s)$. It is readily verifiable
that $|DF(r(s))r'(s)|=\lambda|r'(s)|$. That is, each material element
of $\gamma_{a}$ stretches by a factor of $\lambda$ as advected by
the flow to time $t=b$. Consequently, such total length of any the
curve also changes by a factor of $\lambda$, i.e. $\ell(\gamma_{b})=\lambda\ell(\gamma_{a})$,
where $\ell$ is the length of the curve.

{For $\lambda=1$, this implies that the final length
$\ell(\gamma_{b})$ is equal to the initial length $\ell(\gamma_{a})$,
and hence the $\lambda$ line exactly preserves its arclength. This
is a highly atypical behavior in a turbulent flow, where a typical material
line will stretch exponentially under advection. Yet through any point
in the domain $U$, there will be precisely two material lines preserving
their arclength between the times $a$ and $b$. Such lines are computable
as trajectories of the vector fields $\eta_{1}^{+}$ and $\eta_{1}^{-}$.}

{For $\lambda\neq1$, a similar statement holds for
the subset $U_{\lambda}\subset U$: Passing through any point in $U_{\lambda}$
are two uniformly stretching material lines that stretch by a factor
$\lambda$ between the time $a$ and $b$.}

{(ii) }{\emph{Existence of closed
$\lambda$-lines}}{: Although $\lambda$-lines fill
the set $U_{\lambda}$, most of them are open curves. As shown in
section \S\ref{sec:results}, however, nested families of closed
$\lambda$-lines do arise in two-dimensional turbulence. Members of
such families corresponding to $\lambda=1$ mark the highest possible
degree of coherence in incompressible flows: both of their arclength
and their enclosed area is preserved under material advection. Outermost
members of closed $\lambda$-line families mark Lagrangian vortex
boundaries, the largest possible closed curves that remain coherent
under advection \citep{bhEddy}}.

{(iii)}\emph{ }Relation to KAM tori: In time-periodically
perturbed, two-dimensional Hamiltonian systems, Kolmogorov--Arnold--Moser
(KAM) curves are closed material lines that are mapped exactly into
themselves by the flow in one time-period \cite{guckenheimer}. These
curves, therefore, preserve both their arclength and their enclosed
area in one time period, acting as archetypical coherent Lagrangian
vortex boundaries. In a general, temporally aperiodic velocity field,
closed material lines are no longer mapped into their original position
for any choice of the advection time. A closed $\lambda$-line with
$\lambda=1$, however generalizes the notion of a KAM curve in a finite-time
aperiodic flow, exhibiting both conservation of arclength and enclosed
area between the initial and the final time. In the time-periodic
case, closed $\lambda$-lines with $\lambda=1$ become indistinguishable
from KAM curves when extracted over a time that is a high enough multiple
of the period{{} \citep{geotheory,hadji2013}.}\\

In light of the above discussion, we seek Lagrangian
coherent vortex boundaries as closed $\lambda$-lines. We refer to
closed $\lambda$-lines as \emph{elliptic Lagrangian
coherent structures} (or elliptic LCSs, for short).
In the case $\lambda=1$, they are referred to as \emph{primary}
elliptic LCSs.

\subsection{{Metric interpretation and cosmological analogy}}

{As pointed out in \citet{bhEddy}, elliptic LCSs
bear a mathematical analogy with structures surrounding black holes
in cosmology. Over the subset $U_{\lambda}$ of the flow domain, the
bilinear form 
\[
g_{\lambda}(v,w)=\langle v,E_{\lambda}w\rangle
\]
defines a Lorentzian metric with signature $(-,+)$. The set $U_{\lambda}$
equipped with this metric is a two-dimensional Lorentzian manifold
or }{\emph{space-time}}{{} \citep{beem1996global}.
Unlike in Euclidean geometry, the distance between two distinct points
of this space-time, as measured by its metric $g_{\lambda}$, can
be negative or zero.}

In the language of Lorentzian geometry, the $\lambda$-lines
defined by \eqref{eq:null} can be interpreted as closed null-geodesics
of the metric $g_{\lambda}$.\cite{bhEddy,karrasch2014automated} In cosmology, such surfaces of null-geodesics
with closed space-like projections are called \emph{photon
spheres}.\cite{Claudel2001,bhEddy} They
are composed of periodic light orbits that encircle black holes.

An elliptic LCS, as any closed null-geodesic of the metric $g_{\lambda}$,
must necessarily encircle at least two singularities of $g_{\lambda}$.
This fact considerably simplifies the automated detection of elliptic
LCSs in spatially complex flow data \cite{karrasch2014automated}.

\section{Results and discussion}

{\label{sec:results} We use the method described
in section \S\ref{sec:theory} to identify coherent Lagrangian vortices
in a direct numerical simulation of two-dimensional forced turbulence. }

\subsection{{Numerical method}}

{Consider the Navier--Stokes equations \begin{subequations}
\label{eq:nse} 
\begin{equation}
\partial_{t}u+u\cdot\nabla u=-\nabla p+\nu\Delta u+f,
\end{equation}
\begin{equation}
\nabla\cdot u=0,
\end{equation}
\begin{equation}
u(x,0)=u_{0}(x),
\end{equation}
\end{subequations} where the velocity field $u(x,t)$ is defined
on the two-dimensional domain $U=[0,2\pi]\times[0,2\pi]$ with doubly
periodic boundary conditions.}

{We use a standard pseudo-spectral method with $512$
modes in each direction and $2/3$ dealiasing to solve the above Navier--Stokes
equation with viscosity $\nu=10^{-5}$. The time integration is carried
out over the interval $t\in[0,50]$ (approximately, three eddy-turn-over
times) by a fourth-order Runge-Kutta method with variable step-size
\citep{ode45}. The initial condition $u_{0}$ is the velocity field
of a decaying turbulent flow. The external force $f$ is random in
phase and band-limited, acting on the wave-numbers $3.5<k<4.5$. The
forcing amplitude is time-dependent balancing the instantaneous enstrophy
dissipation $\nu\int k^{2}Z(k,t)\,\id k$ where $Z(k,t):=\frac{1}{2}\int_{|\vc k|=k}|\hat{\omega}(\vc k,t)|^{2}\,\id S(\vc k)$
with $\hat{\omega}(\cdot,t)$ being the Fourier transform of the instantaneous
vorticity $\omega(\cdot,t)=\nabla\times u(\cdot,t)$.}

{In two dimensions, the energy injected into the
system by the forcing is mostly transferred to larger scales through
a nonlinear process \citep{kraichnan1967,merilees}. In order to prevent
the energy accumulation at largest available scales over time, a linear
damping is usually added to the Navier--Stokes equation to dissipate
the energy at large scales \citep{bofetta-ekman,tsang2005}. However,
for the time scales considered here, the energy accumulation is not
an issue and hence the linear damping will be omitted.}

{The theory reviewed in Section \S\ref{sec:theory}
does not assume a particular governing equation for the velocity field
$u(x,t)$. Thus, it can be applied to any two-dimensional velocity
field, given as numerical solution of a partial differential equation
or by direct measurements. In particular, it can be applied to Lagrangian
vortex detection for the solutions of the Navier--Stokes equation
\eqref{eq:nse}. To detect the Lagrangian vortex boundaries, we take
the following steps: }
\begin{enumerate}
\item {Solve the Navier--Stokes equation \eqref{eq:nse}
as discussed above to get the velocity field $u(x,t)$ over the time
interval $t\in[0,50]$ and a uniform $512\times512$ spatial grid
over the domain $x\in U=[0,2\pi]\times[0,2\pi]$. The temporal resolution
of the velocity field is $251$ such that two consecutive time slices
are $\Delta t=0.2$ apart.}
\item {Advect each grid point according to the differential
equation \eqref{eq:dynsys} from time $t=0$ to time $t=50$ to construct
the flow map $F$ such that $F(x_{a})=x_{b}$ for any grid point $x_{a}$. }
\item {Construct an approximation of the deformation gradient
$DF$ by finite differences. To increase the finite difference accuracy,
we use the auxiliary grid method introduced in \citet{computeVariLCS}.
The chosen auxiliary grid distance is $10^{-3}$. }
\item {Construct the right Cauchy--Green strain tensor
$C(x_{a})=[DF(x_{a})]^{\top}DF(x_{a})$ for each grid point $x_{a}$.
Compute the eigenvalues $\{\lambda_{1},\lambda_{2}\}$ and the corresponding
eigenvectors $\{\xi_{1},\xi_{2}\}$ of $C(x_{a}).$ }
\item {Seek closed orbits of the one-parameter families
of vector fields $\eta_{\lambda}^{\pm}$ defined in \eqref{eq:eta}.
For detecting these closed orbits, we use the automated algorithm
developed in \citet{bhEddy}. }
\end{enumerate}
{We detect the Lagrangian vortex boundaries as outermost
elliptic LCSs, i.e., maximal limit cycles of $\eta_{\lambda}^{\pm}$.
In the following, we present a detailed analysis of these vortex boundaries
and compare them to those suggested by alternative Eulerian and Lagrangian
indicators.}

\subsection{{Lagrangian coherent vortex analysis}}

\begin{figure}
\centering
\subfloat[]{\includegraphics[width=0.47\textwidth]{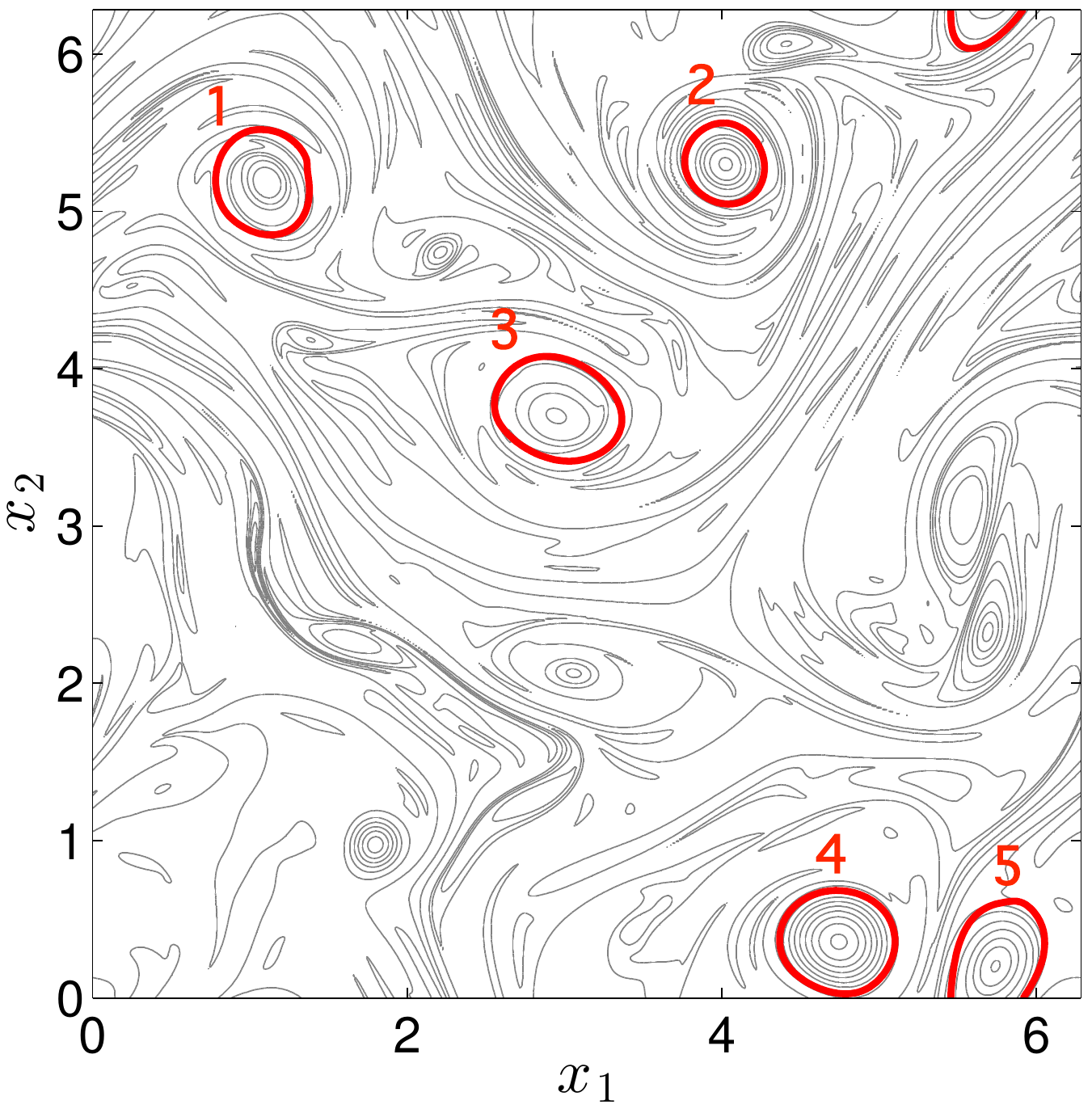}
\includegraphics[width=0.47\textwidth]{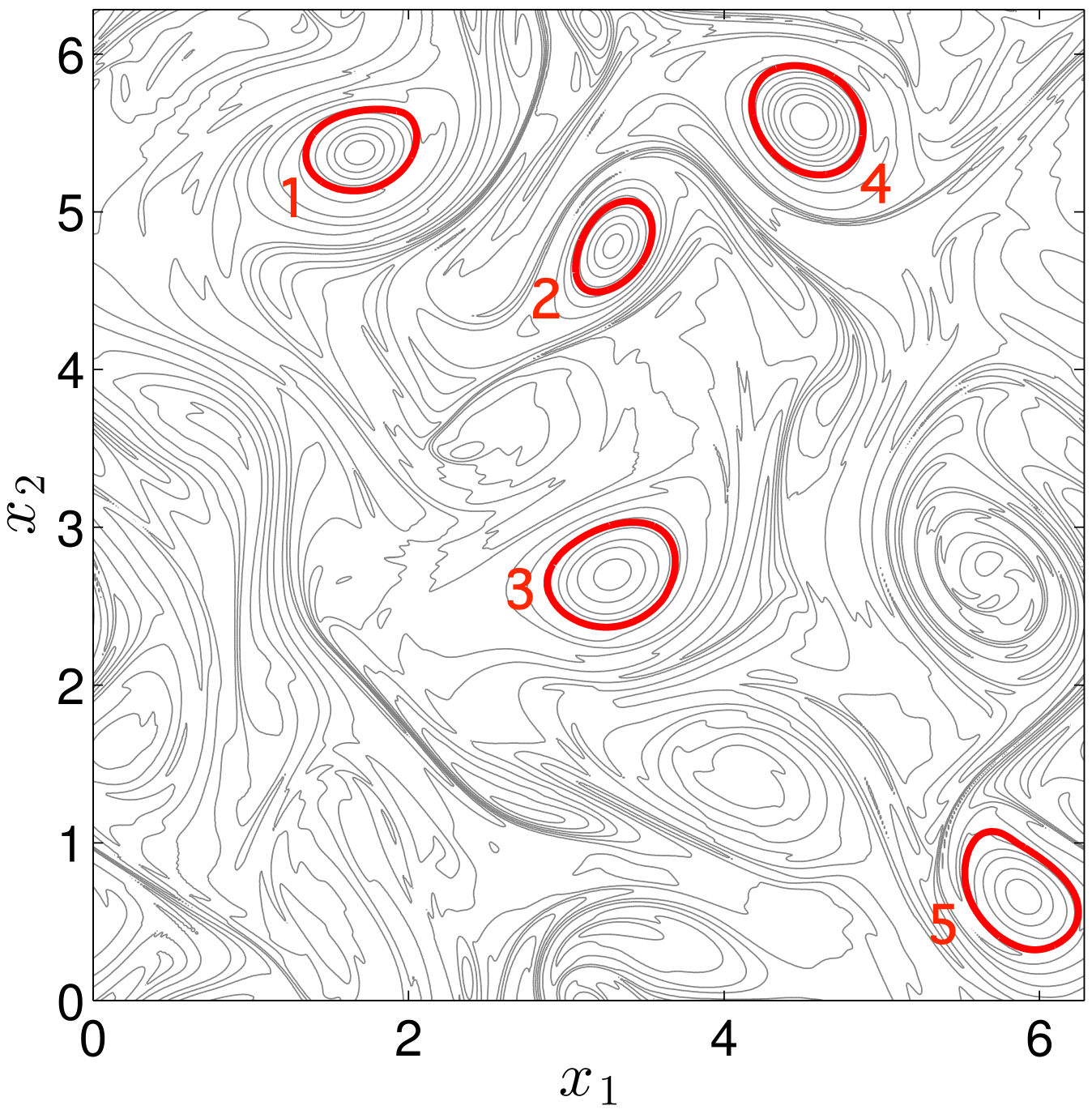}}\\
\subfloat[]{\includegraphics[width=0.47\textwidth]{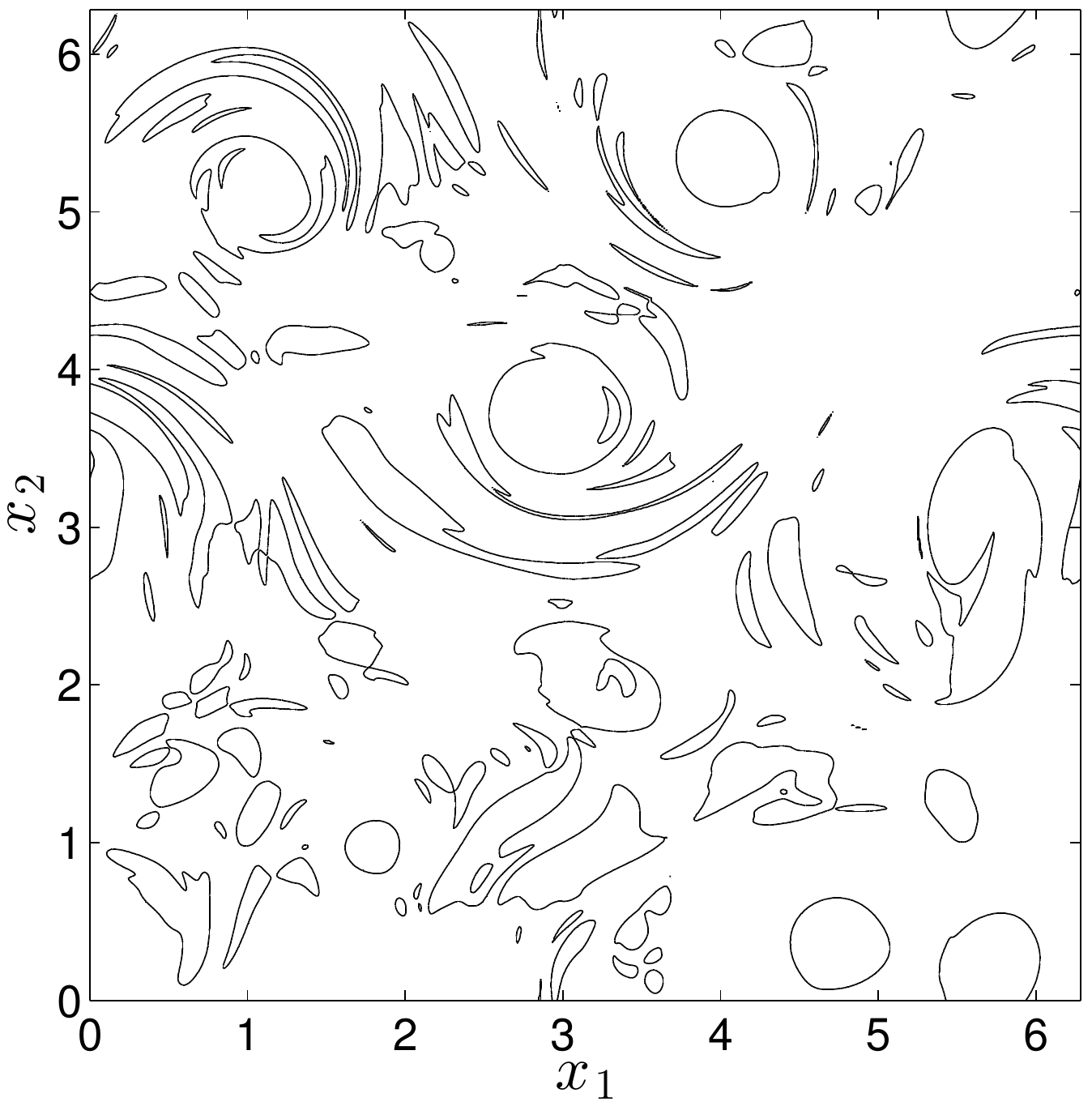}
\includegraphics[width=.47\textwidth]{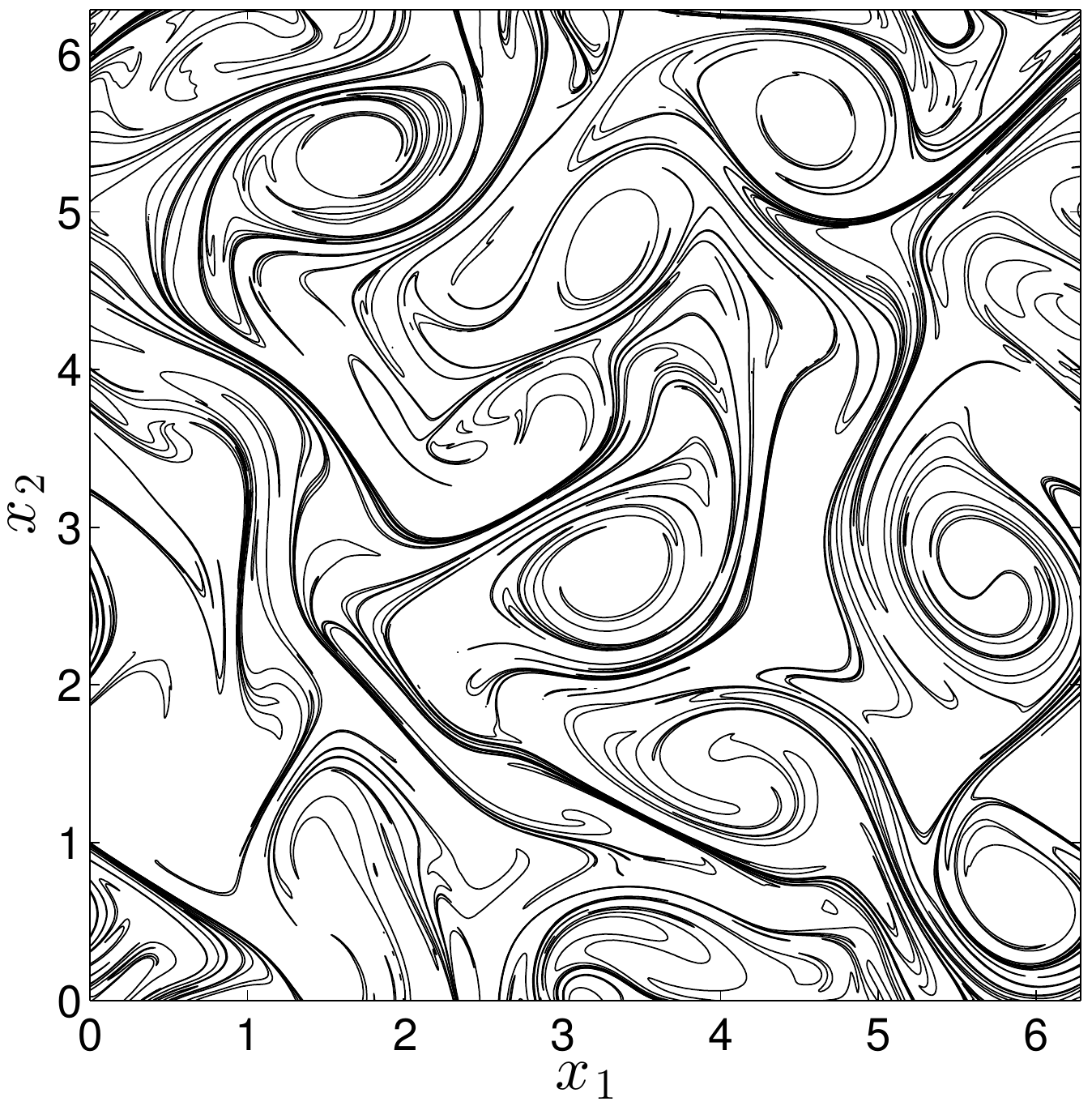}}
\caption{(a) Lagrangian vortex boundaries (red) at time $t=0$ (left) and $t=50$
(right). The vorticity contours are shown in gray in the background.
The vorticity contours are distributed as $-1:0.1:1$ at time $t=0$
and as $-1.5:0.15:1.3$ at time $t=50$. The coherent vortices are
numbered in order to facilitate their identification at the two time-instances.
(multimedia view) (b) Zero level curves of the Okubo-Weiss parameter
at $t=0$ (left) and their advected positions at time $t=50$ (right).}
\label{fig:nostrlines+vorticity}
\end{figure}

{Figure \ref{fig:nostrlines+vorticity}a (left) shows
the boundaries (red) of Lagrangian coherent vortices superimposed
on the contours of the Eulerian vorticity (gray) at time $t=0$. The
boundaries are found as the outermost elliptic LCSs, i.e., maximal
limit cycles of the vector fields $\eta_{\lambda}^{\pm}$ (see Eq.
(\ref{eq:eta})). The advected coherent vortex boundaries at time
$t=50$ are shown in figure \ref{fig:nostrlines+vorticity}a (right)
along with the corresponding instantaneous vorticity field. By construction,
these Lagrangian vortex boundaries resist straining and filamentation
under advection {[}see Fig. \ref{fig:nostrlines+vorticity} (multimedia
view){]}. In the following, the vortex numbers refer to the numbering
in figure \ref{fig:nostrlines+vorticity}a.}

{For basic reference, we also plot the zero level-curves
of the most often used two-dimensional vortex diagnostic, the Okubo--Weiss
parameter\citep{okubo,weiss_okubo}, at $t=0$ (Figure \ref{fig:nostrlines+vorticity}b,
left). We then advect these contours to $t=50$ (Figure \ref{fig:nostrlines+vorticity}b,
right). Clearly, the Okubo--Weiss zero curves deform significantly,
and hence cannot be considered as approximations to coherent material
vortex boundaries. This is in line with similar observations made
by earlier studies (see, e.g., \citet{pasquero2001}, \citet{isern2006},
\citet{Henson2008} and \citet{javier_agulhas}). We present the definition
and a detailed analysis of the Okubo--Weiss parameter in Section \ref{sec:compare}. }

\begin{figure}[t!]
\begin{center}
\subfloat[]{{\includegraphics[width=0.45\textwidth]{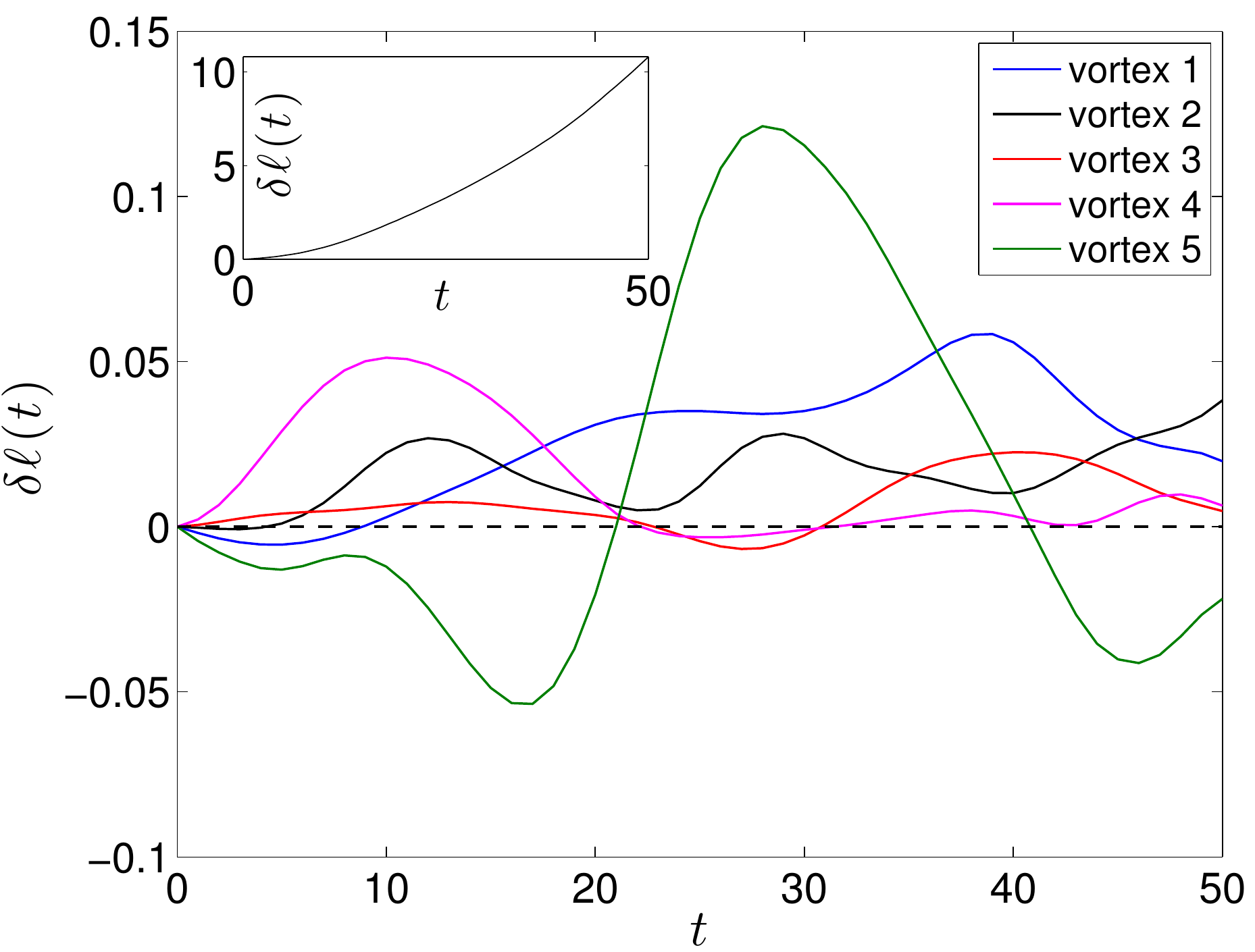}}}
\hspace{.03\textwidth}
\subfloat[]{{\includegraphics[width=0.5\textwidth]{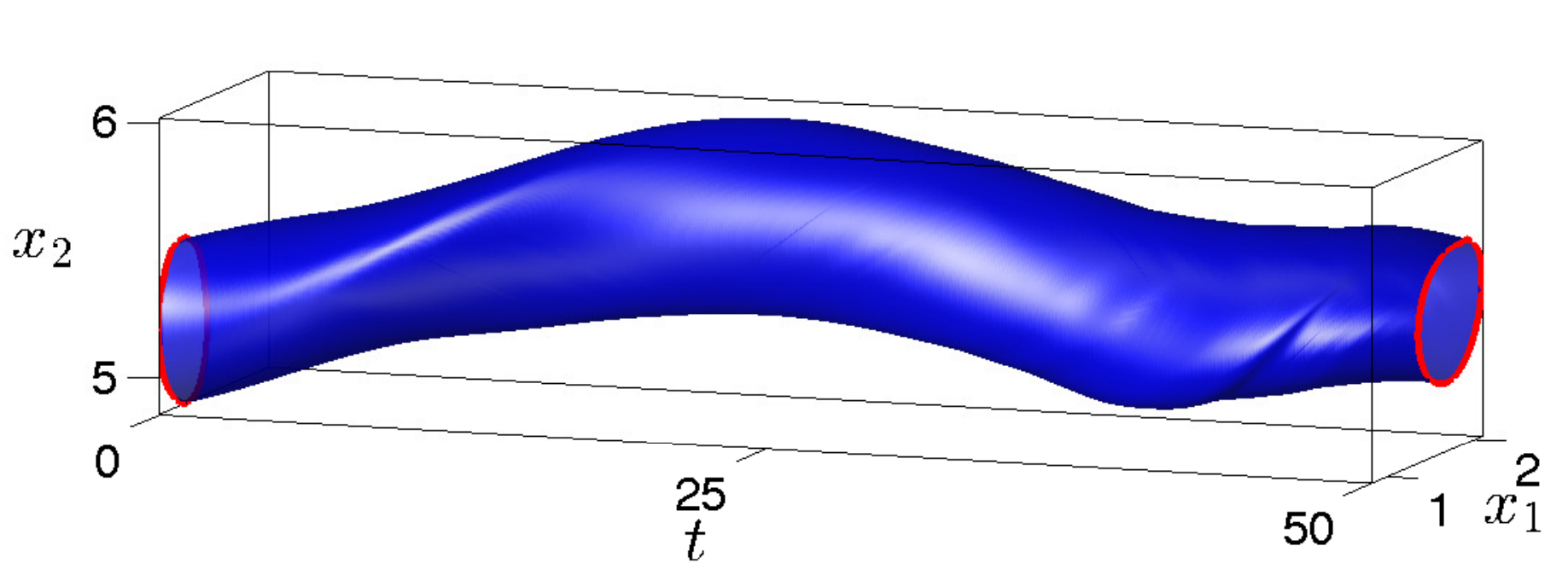}}}
\end{center}
\caption{(a) The relative deformation $\delta\ell$ as a function of time for the primary
elliptic LCSs. The inset shows the relative stretching for a typical
closed material line over the same time interval. (b) The Lagrangian
vortex 1 in the extended phase space. The tube is created from the
advection of the vortex boundary under the flow.}
\label{fig:length_of_cshrlines}
\end{figure}

{Returning to the analysis of the elliptic LCSs,
figure \ref{fig:length_of_cshrlines} shows the relative stretching
$\delta\ell(t):=(\ell(t)-\ell(a))/\ell(a)$ of the primary elliptic
LCSs over the time interval $t\in[0,50]$. Here, $\ell(t)$ denotes
the length of a material line at time $t$. In principle, the initial
and the final lengths of a primary elliptic LCS must be exactly equal,
resulting in zero relative stretching at time $t=50$. In practice,
a deviation of at most $4\%$ is observed from this ideal limit owing
to numerical errors. The inset of figure \ref{fig:length_of_cshrlines}
shows the relative stretching of a typical non-coherent iso-vorticity
line. Unlike the coherent vortices, the relative stretching for a
general material curve increases exponentially, with its final value
at least an order of magnitude larger than that for a coherent vortex.}

\begin{figure}[t!]
\centering
\subfloat[]{\includegraphics[width=0.45\textwidth]{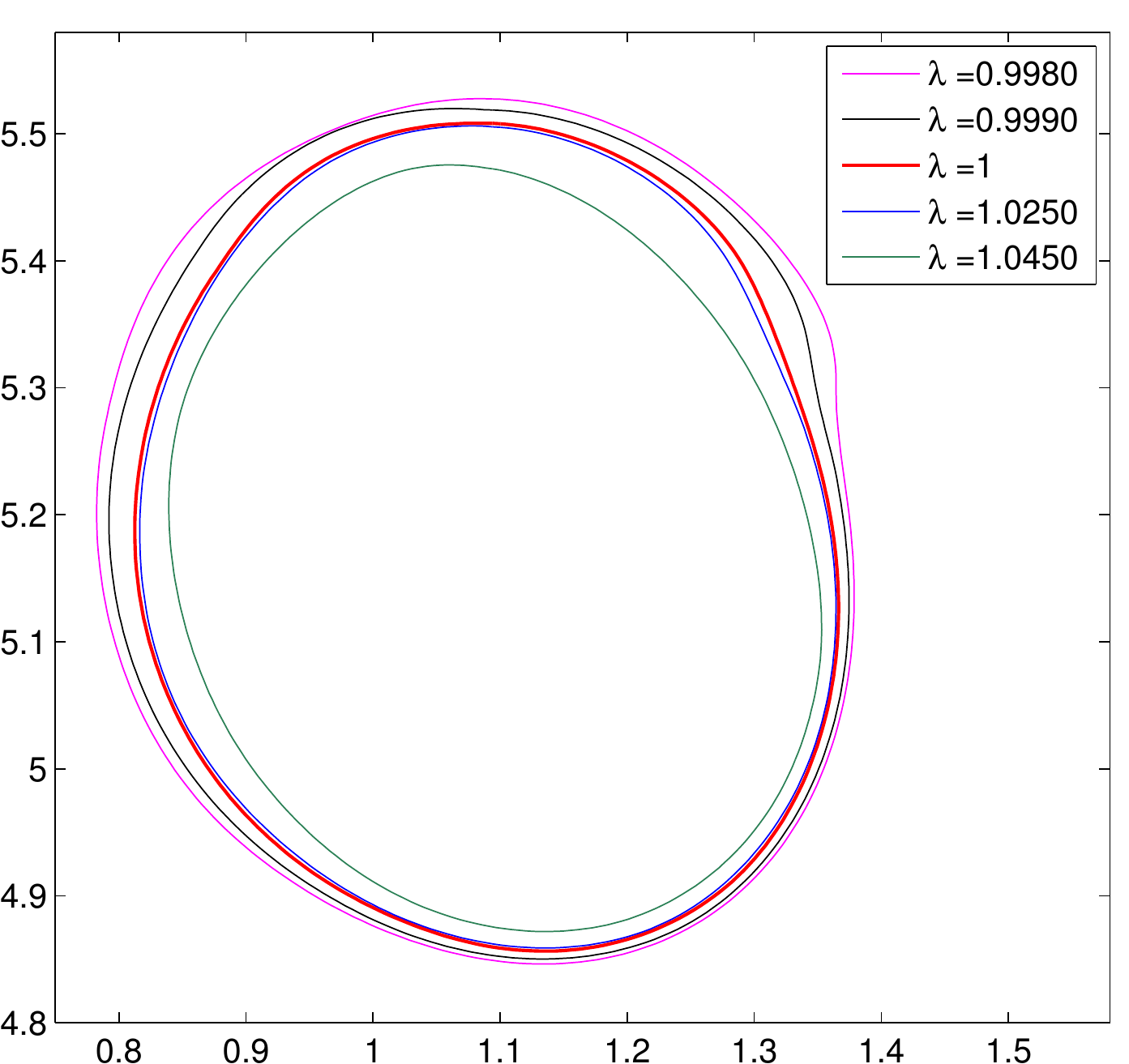}}
\hspace{.08\textwidth}
\subfloat[]{\includegraphics[width=0.45\textwidth]{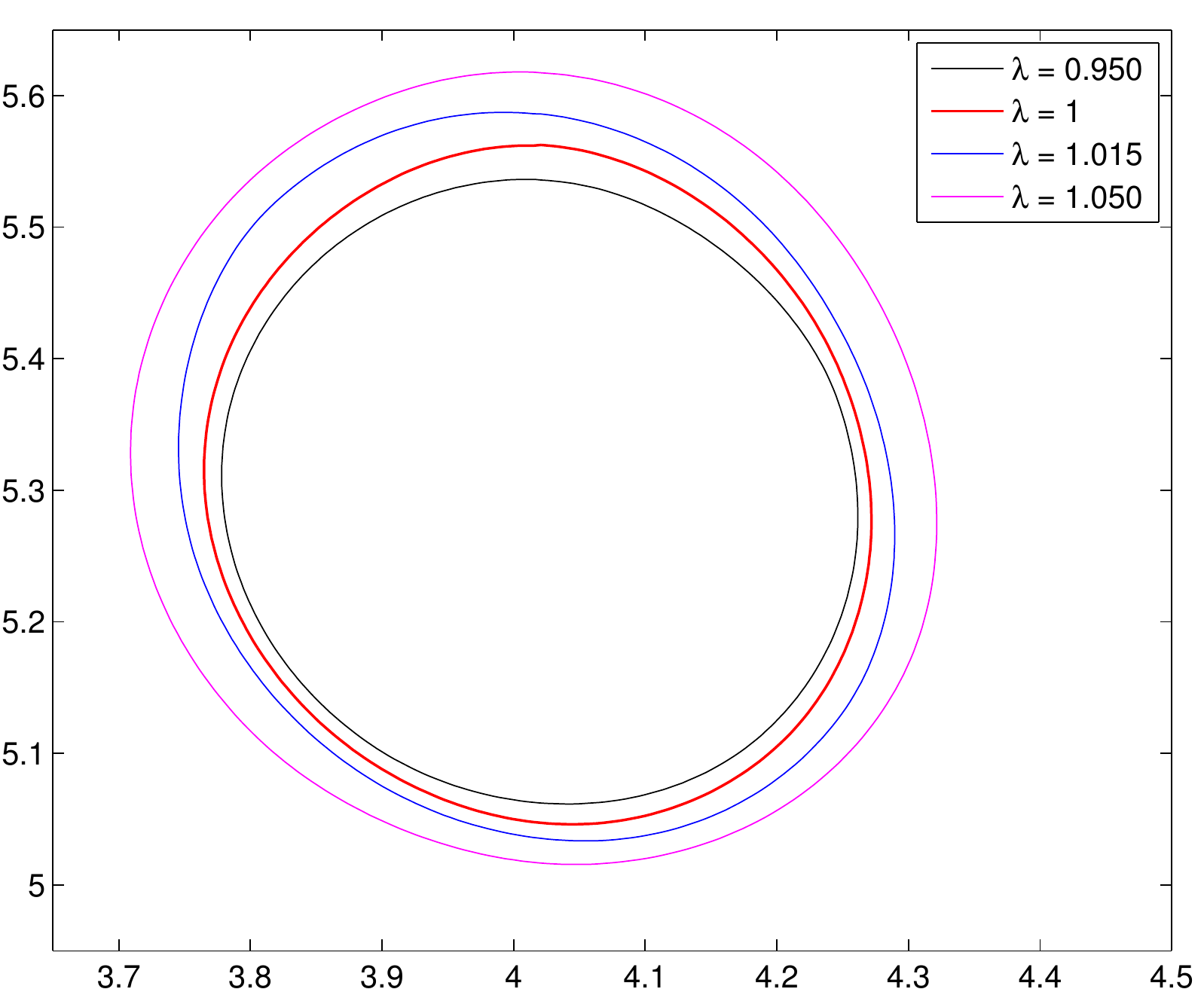}}
\caption{Elliptic LCSs (i.e., closed $\lambda$-lines) around vortex 1 (a)
and vortex 2 (b).}
\label{fig:lamStrlines_nested}
\end{figure}

{As mentioned in section \S\ref{sec:theory}, coherent
material vortex boundaries are formed by a nested set of elliptic
LCSs (i.e., closed $\lambda$-lines). Figure \ref{fig:lamStrlines_nested}
shows two of the coherent vortices and their corresponding $\lambda$-lines.
We find that for vortex 1, the secondary elliptic LCSs with $\lambda>1$
lie in the interior of the primary elliptic LCS (i.e., the closed
$\lambda$-line with $\lambda=1$). For all other coherent vortices
of figure \ref{fig:nostrlines+vorticity}, the secondary elliptic
LCSs with $\lambda>1$ lie in the exterior of the primary elliptic
LCS. In all five cases, values of $\lambda$ for which an elliptic
LCS exists are close to 1, ranging in the interval $0.94\leq\lambda\leq1.05$.}

The majority of vortices appearing in figure \ref{fig:nostrlines+vorticity}
are not coherent in the Lagrangian frame, and hence no elliptic LCSs
were found around them. Some of the non-coherent vortices are trapped
in a hyperbolic region, experiencing substantial straining over time.
Others undergo a merger process where a larger vortex is created from
two smaller co-rotating vortices. Each smaller vortex deforms substantially
during the merger. The merged vortex may or may not remain coherent
for later times.

\begin{figure}[t!]
\centering
\subfloat[]{\includegraphics[width=0.4\textwidth]{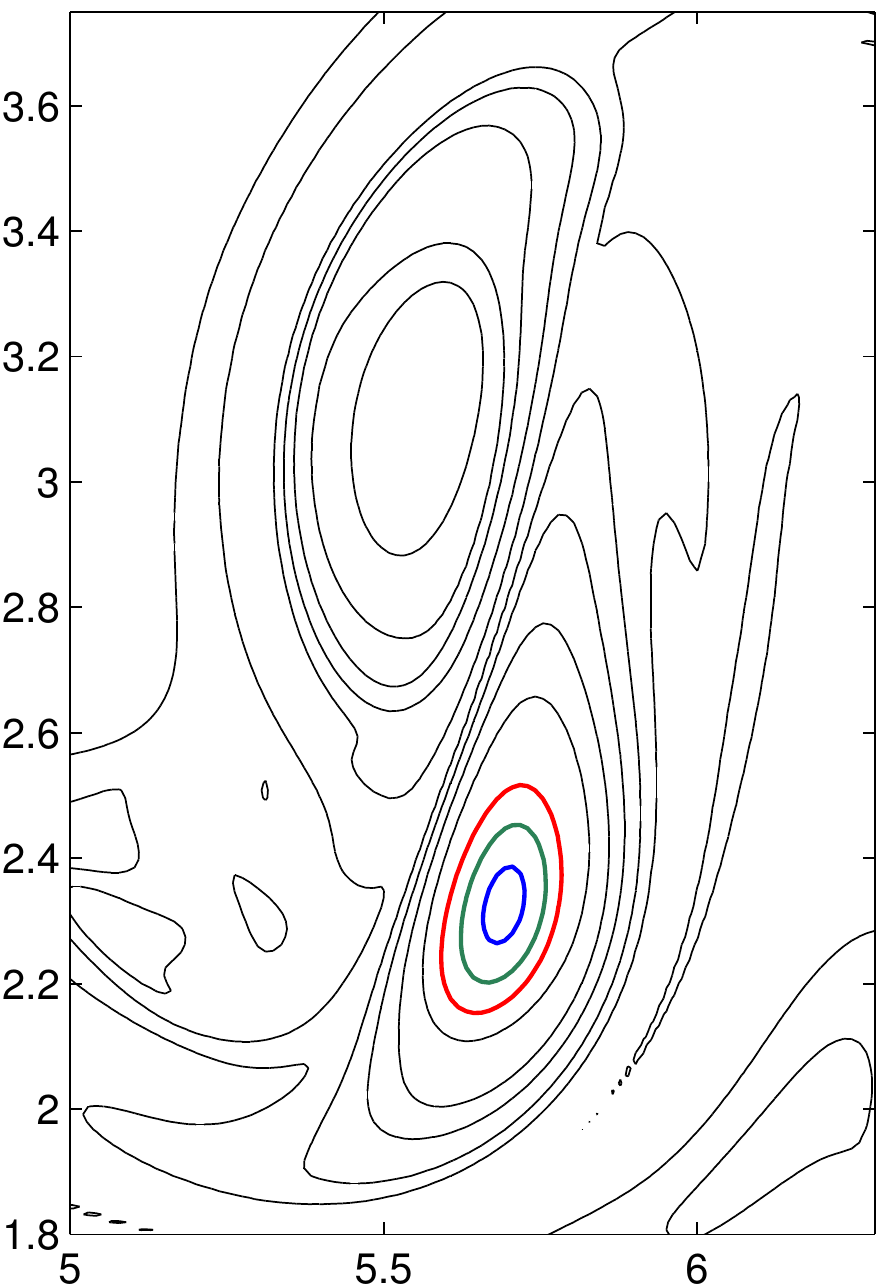}}
\hspace{.1\textwidth}
\subfloat[]{\includegraphics[width=0.45\textwidth]{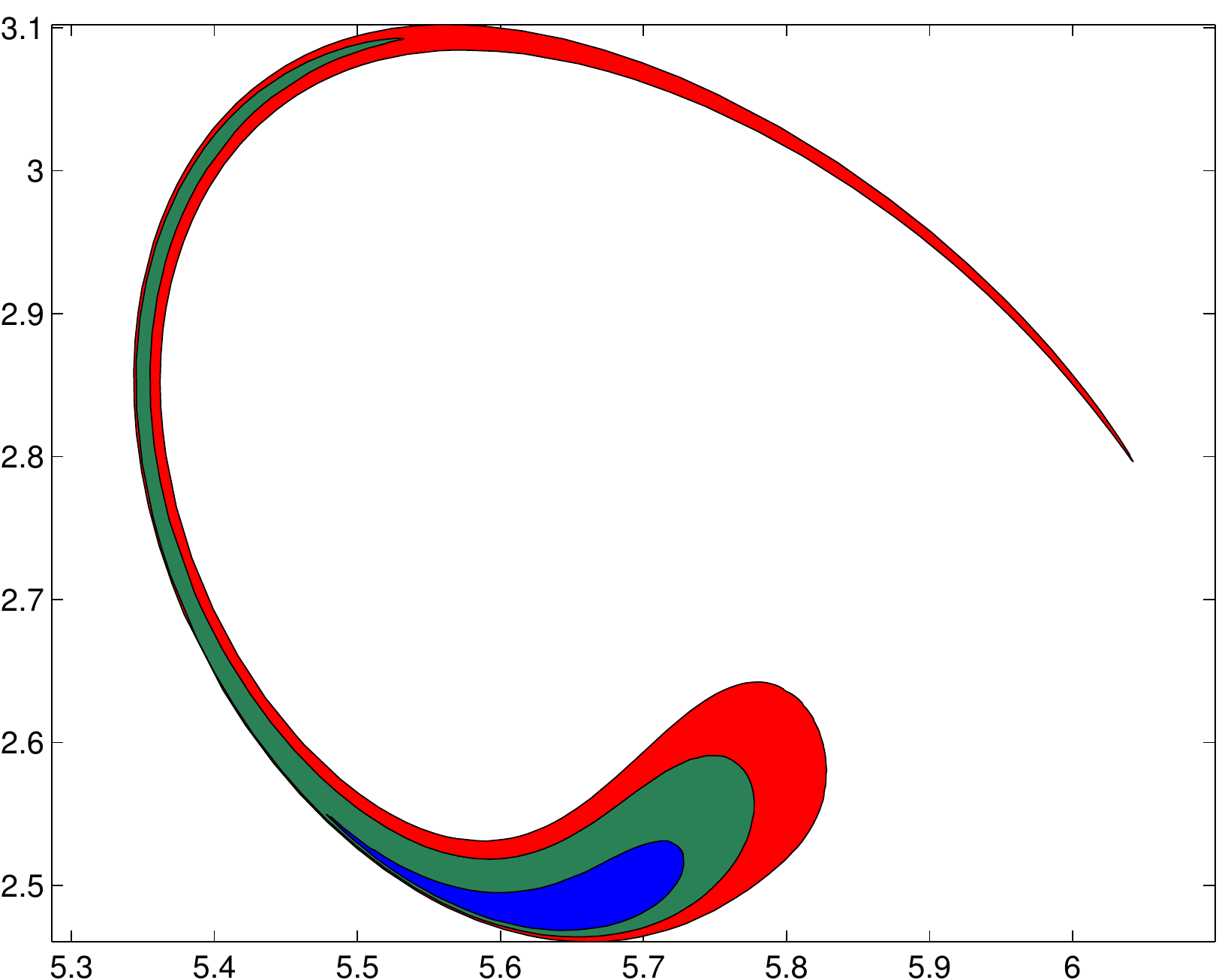}}
\caption{(a) Vortex contours at $t=0$ for two non-coherent vortices that merge
as one later in time. To demonstrate the deformation of the vortices
we monitor the advection of three vorticity contours. The contour
values are $0.6$ (red), $0.7$ (green) and $0.8$ (blue). (b) The
selected contours advected to time $t=50$ and filled with their corresponding
color.}
\label{fig:nonCoherentVort_01}
\end{figure}

\begin{figure}[t!]
\centering\subfloat[]{\includegraphics[width=0.55\textwidth]{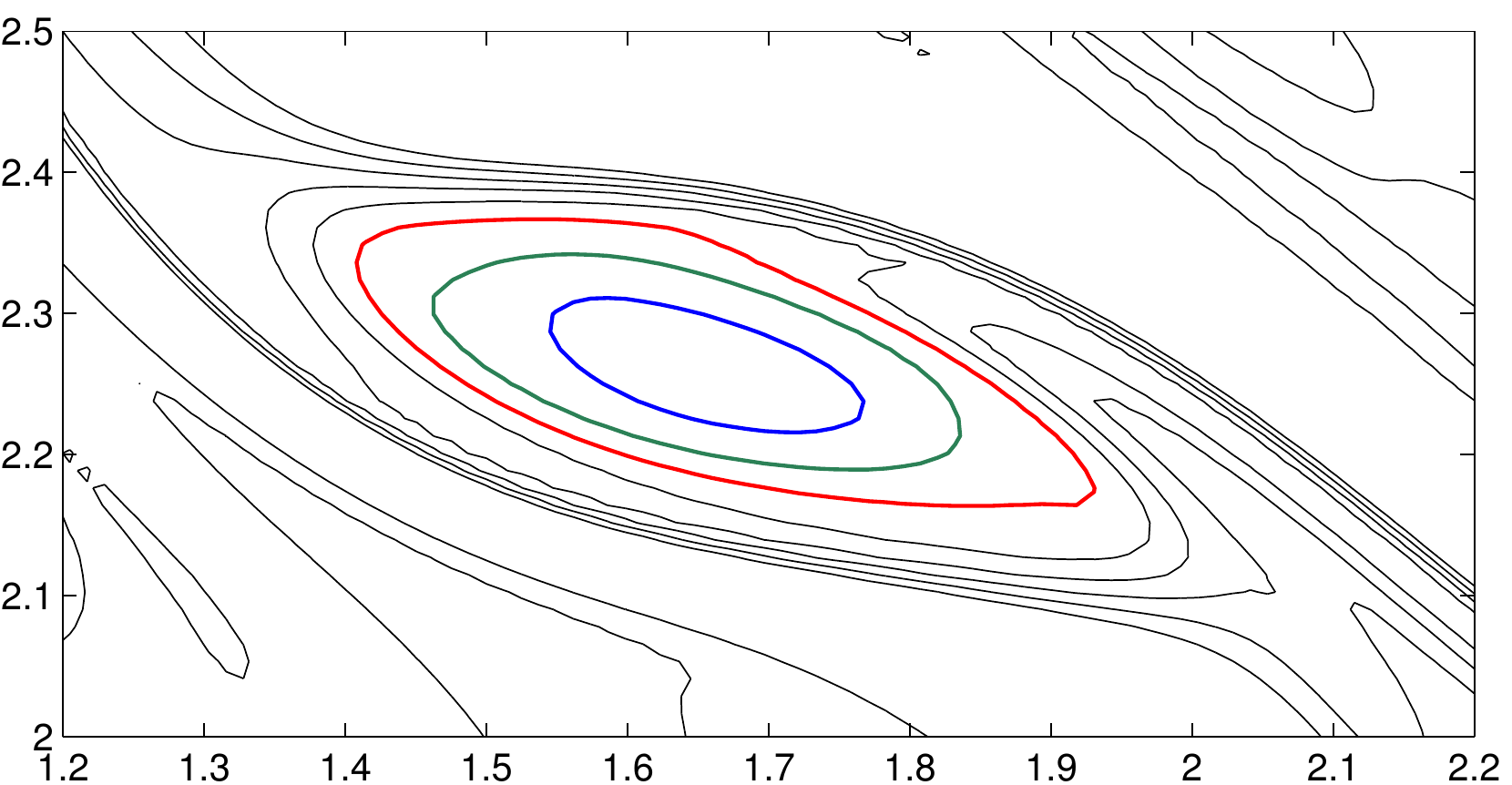}}
\hspace{.04\textwidth}
\subfloat[]{\includegraphics[width=0.4\textwidth]{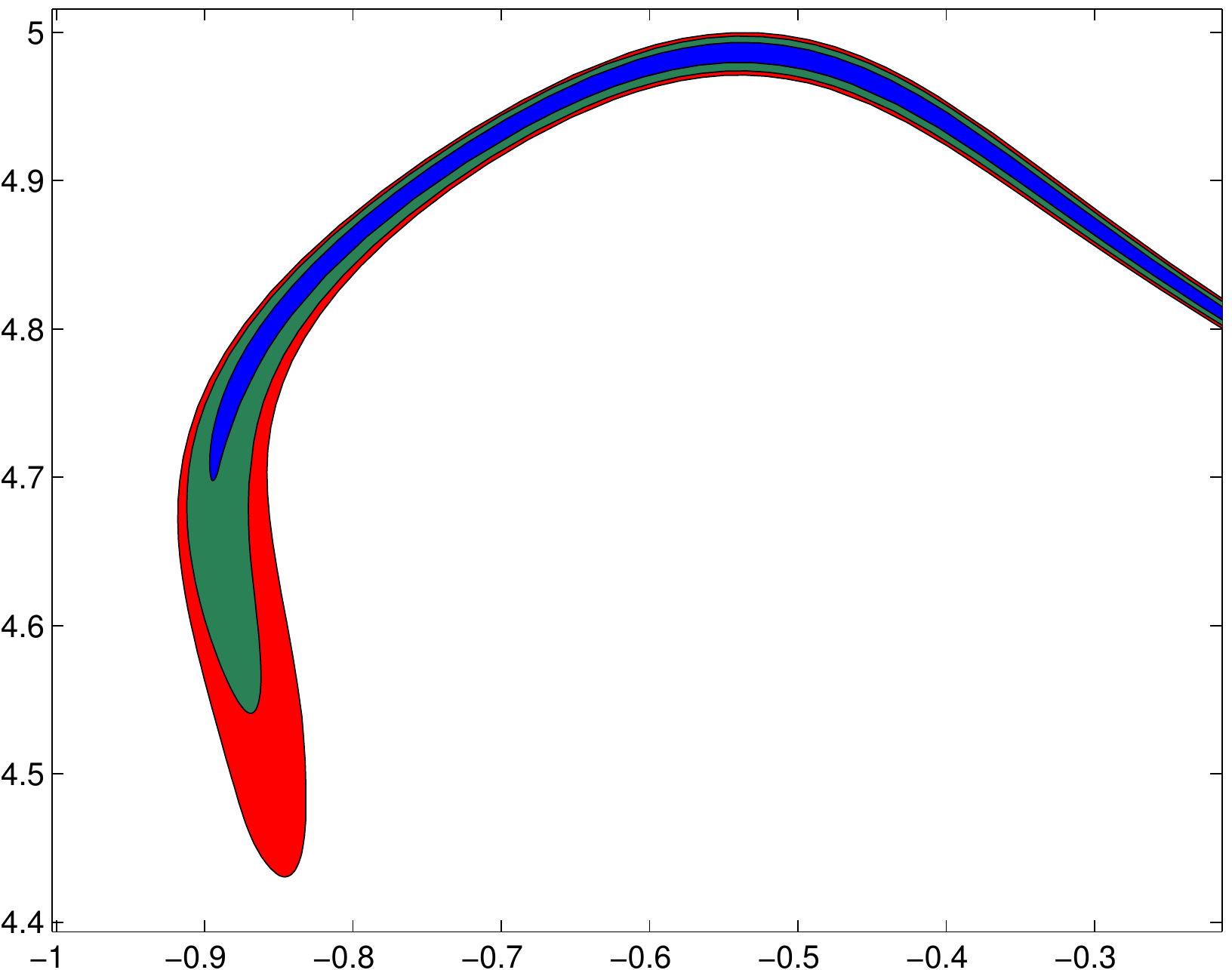}}
\caption{(a) Vortex contours at $t=0$ for a non-coherent vortex trapped in
a straining field. The contours of vorticity with values $0.25$ (red),
$0.3$ (green) and $0.35$ (blue) are marked. (b) The selected contours
advected to time $t=50$ and filled with their corresponding color.
Only part of the advected image is shown.}
\label{fig:nonCoherentVort_02}
\end{figure}

Figure \ref{fig:nonCoherentVort_01} focuses on
one Eulerian vortex undergoing a merger process. To illustrate the
deformation of this vortex, we take three vorticity contours at time
$t=0$ near the center of the vortex. Selected vorticity contours
are then advected to the final time $t=50$, showing the resulting
deformation of the vortex core. Figure \ref{fig:nonCoherentVort_02}
shows a similar analysis for a non-coherent vortex trapped in a uniformly
hyperbolic region of the flow. Hyperbolicity produces stretching of
vorticity gradients resulting in smearing of the vortex.
\begin{figure}[t!]
\centering \includegraphics[width=0.65\textwidth]{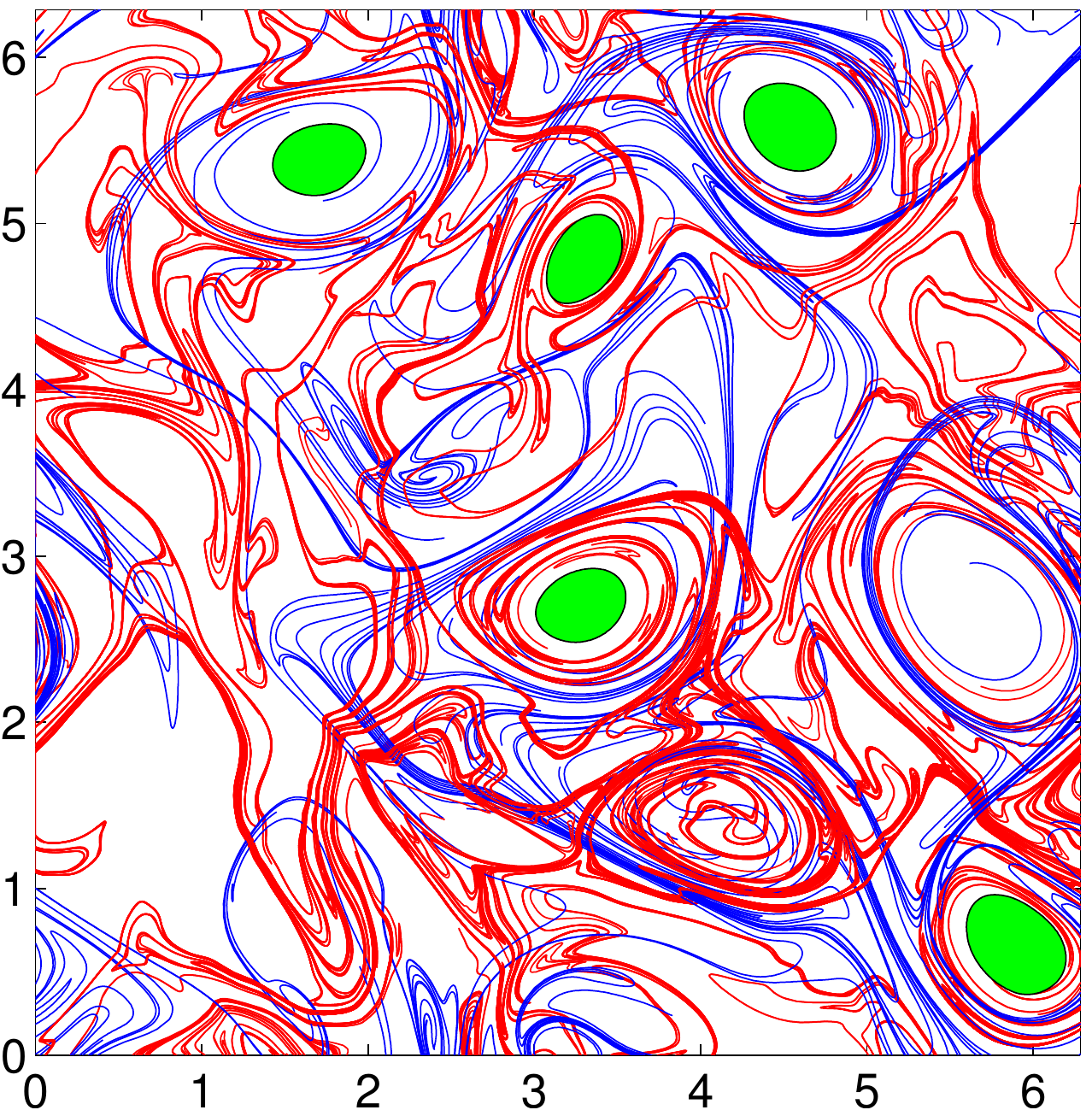}
\caption{Generalized stable (red) and unstable (blue) manifolds. The coherent
Lagrangian vortices (green), i.e. generalized KAM regions, are not
penetrated by these manifolds. The manifolds and the KAM regions are
shown at $t=50$.}
\label{fig:hypLCS}
\end{figure}

Figure \ref{fig:hypLCS} shows the generalized stable
and unstable manifolds obtained by the geodesic theory of Lagrangian
coherent structures \citep{geotheory,shearless}, using the computational
method described in \citet{stretchlines}. These stable and unstable
manifolds are, respectively, the most repelling and attracting material
lines that form the skeleton of turbulent mixing patterns. The exponential
attraction and repulsion generated by these manifolds leads to smearing
of most fluid regions that appear as vortices in instantaneous streamline
and vorticity plots. By contrast, the coherent Lagrangian vortices
we identify remain immune to straining.

\subsection{{Optimality of coherent vortex boundaries}}

{Here we examine the optimality of vortex boundaries
obtained as outermost elliptic LCSs. The optimal boundary of a coherent
vortex can be defined as a closed material line that encircles the
largest possible area around the vortex and shows no filamentation
over the observational time period. We seek to illustrate that outermost
elliptic LCSs mark such optimal boundaries.}

\begin{figure}[t!]
\centering
\subfloat[]{\includegraphics[width=0.4\textwidth]{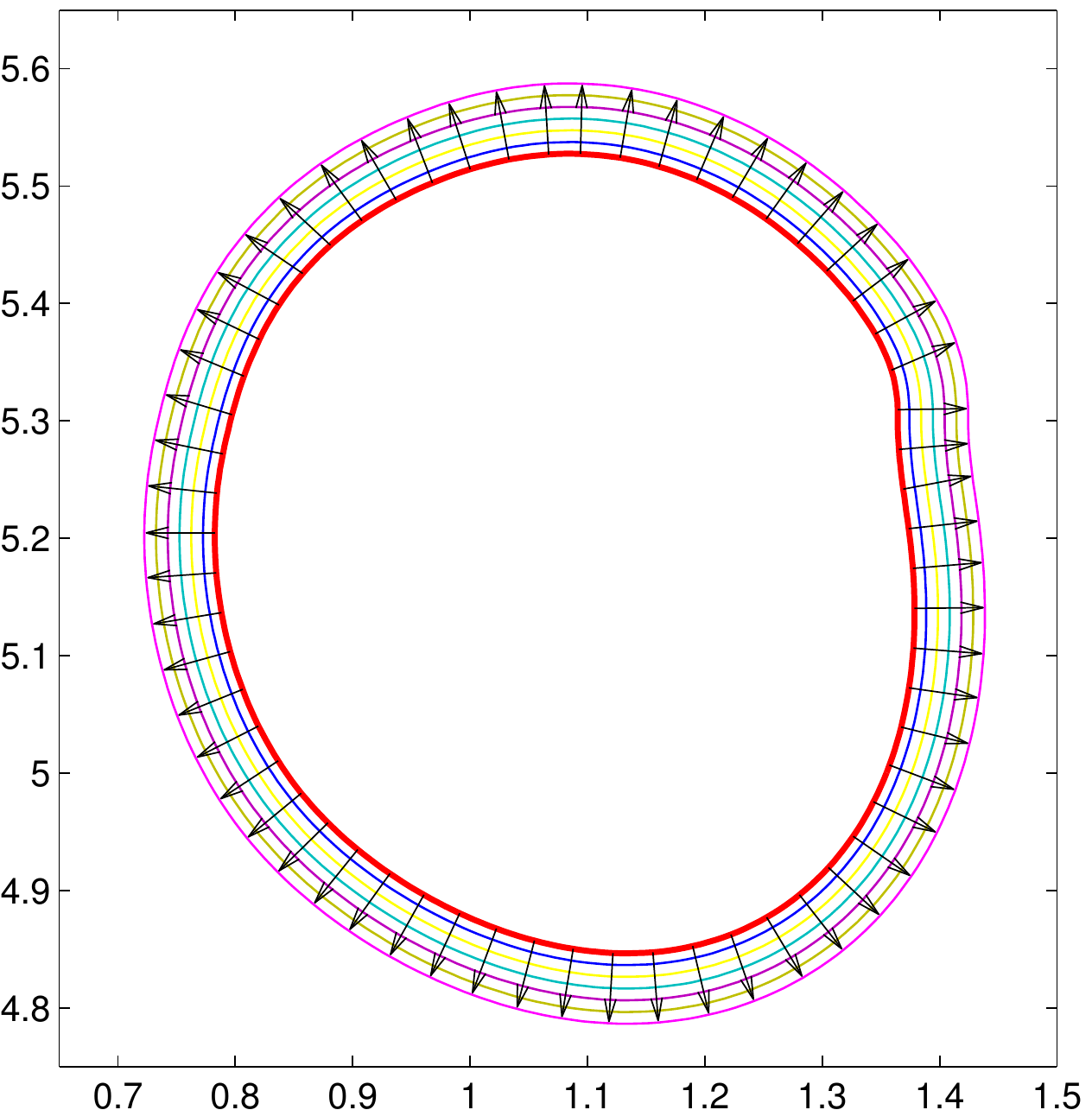}}
\hspace{.04\textwidth}
\subfloat[]{\includegraphics[width=0.55\textwidth]{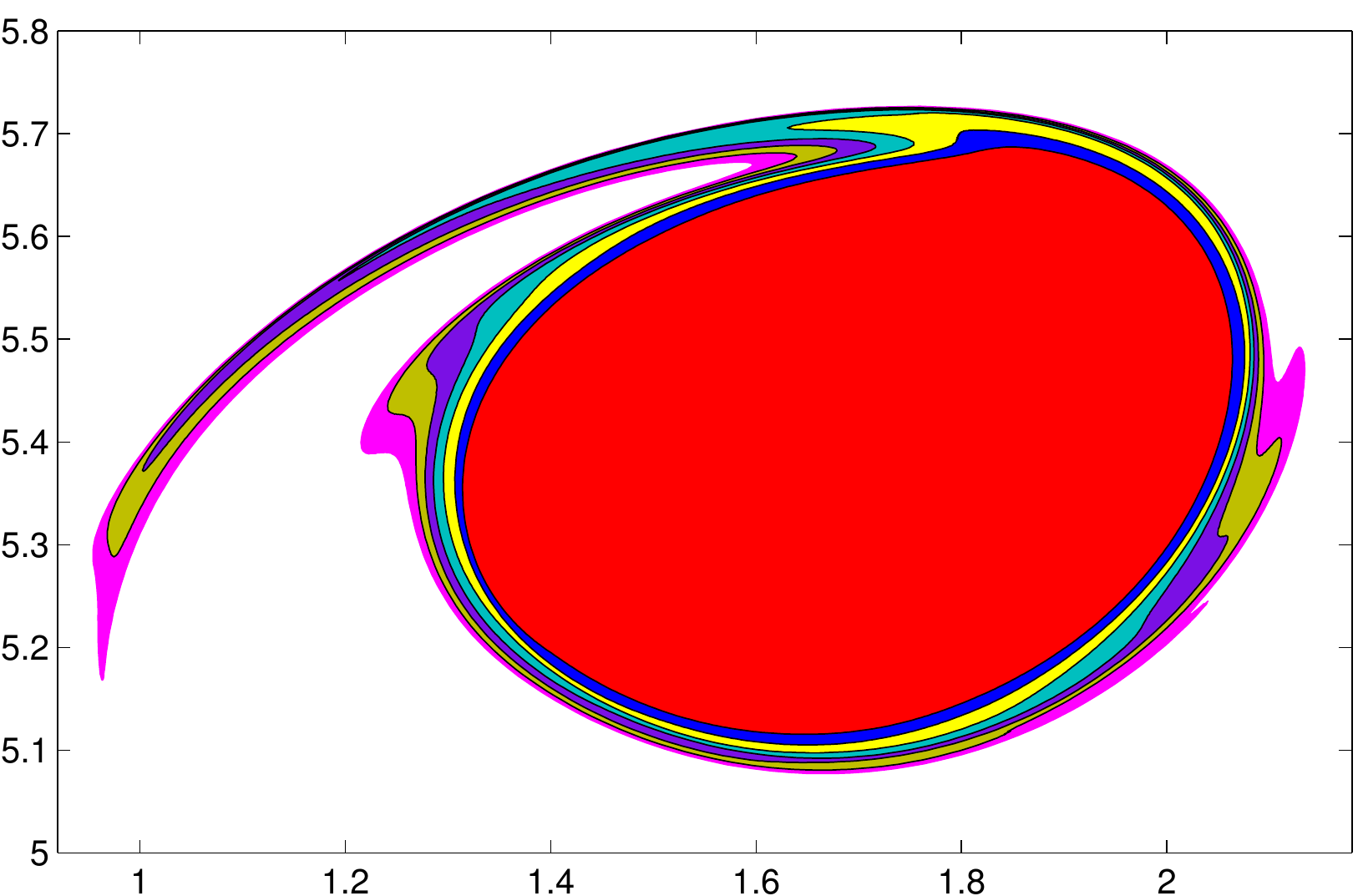}}
\caption{(a) The outermost elliptic LCS (red) encircling vortex 1 of figure
\ref{fig:nostrlines+vorticity} and its outer normal perturbations.
The perturbation parameter ranges between $0.01$ and $0.06$. (b)
The advected image of the elliptic LCS and its normal perturbations
at time $t=50$. Each advected image is filled with its corresponding
color from panel (a)}
\label{fig:optimality}
\end{figure}

{To this end, we consider a class of perturbations
to the outermost elliptic LCS of vortex 1 corresponding to $\lambda=0.998$.
The perturbations are in the direction of the outer normal of the
elliptic LCS. The amount of perturbation ranges between $0.01$ and
$0.06$ (i.e., $1.5\%$ to $10\%$ of the diameter of the elliptic
LCS). We then advect the vortex boundary and its perturbations to
the final time $t=50$ (see figure \ref{fig:optimality}b). The perturbed
curves visibly depart from the coherent core marked by the red elliptic
LCS. Our findings are similar for all other coherent vortices (not
shown here). }

\subsection{Comparison with Eulerian and Lagrangian vortex indicators}\label{sec:compare}

Several diagnostics have been
previously proposed to identify vortex boundaries. Among the Eulerian
indicators are the vorticity criterion of \citet{mcwilliams_vort},
the Okubo-Weiss (OW) criterion \cite{okubo,weiss_okubo} and the
modified OW criterion of \citet{hua1998}. These Eulerian methods
are non-objective (frame-dependent), instantaneous in nature, and
are generally used in practice with tunable thresholds. For all these
reasons, they have little chance to capture long-term coherence in
the Lagrangian frame. Nevertheless, they are broadly believed to be
good first-order indicators of coherence in the flow.

\begin{figure}[t!]
\centering
\subfloat[]{{\includegraphics[width=0.4\textwidth]{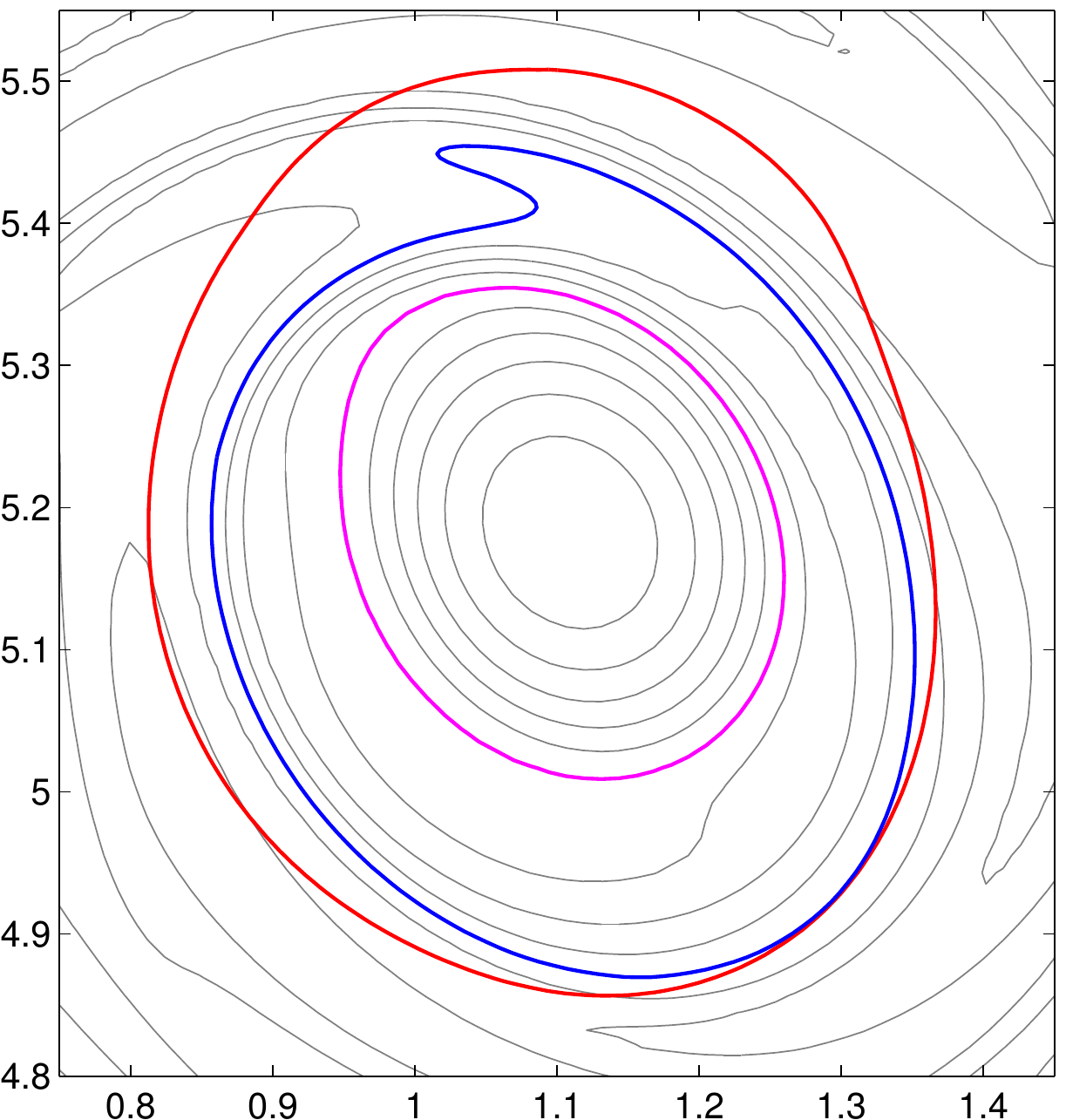}\hspace{.05\textwidth}
\includegraphics[width=0.42\textwidth]{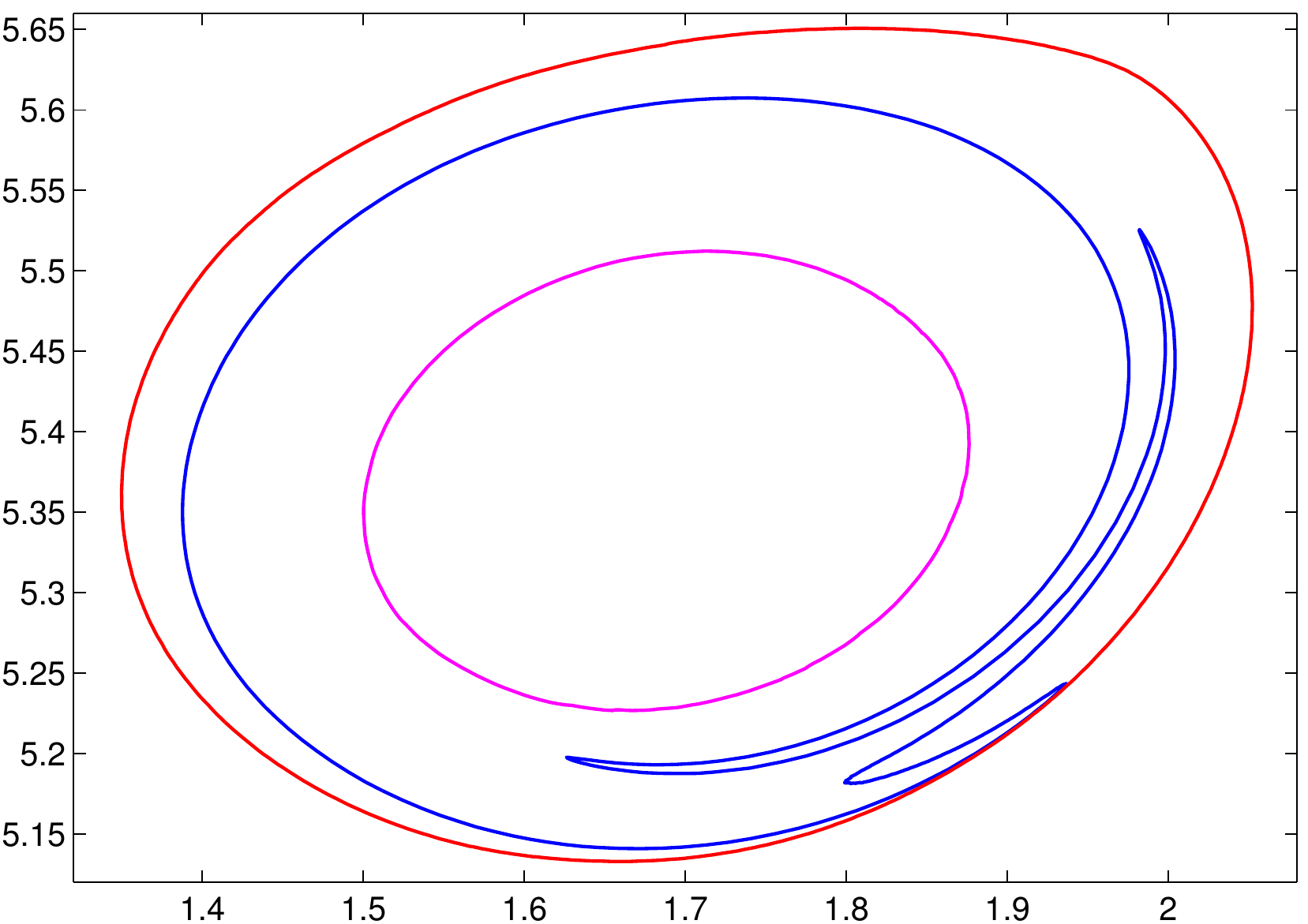}}}\\
\subfloat[]{\includegraphics[width=0.45\textwidth]{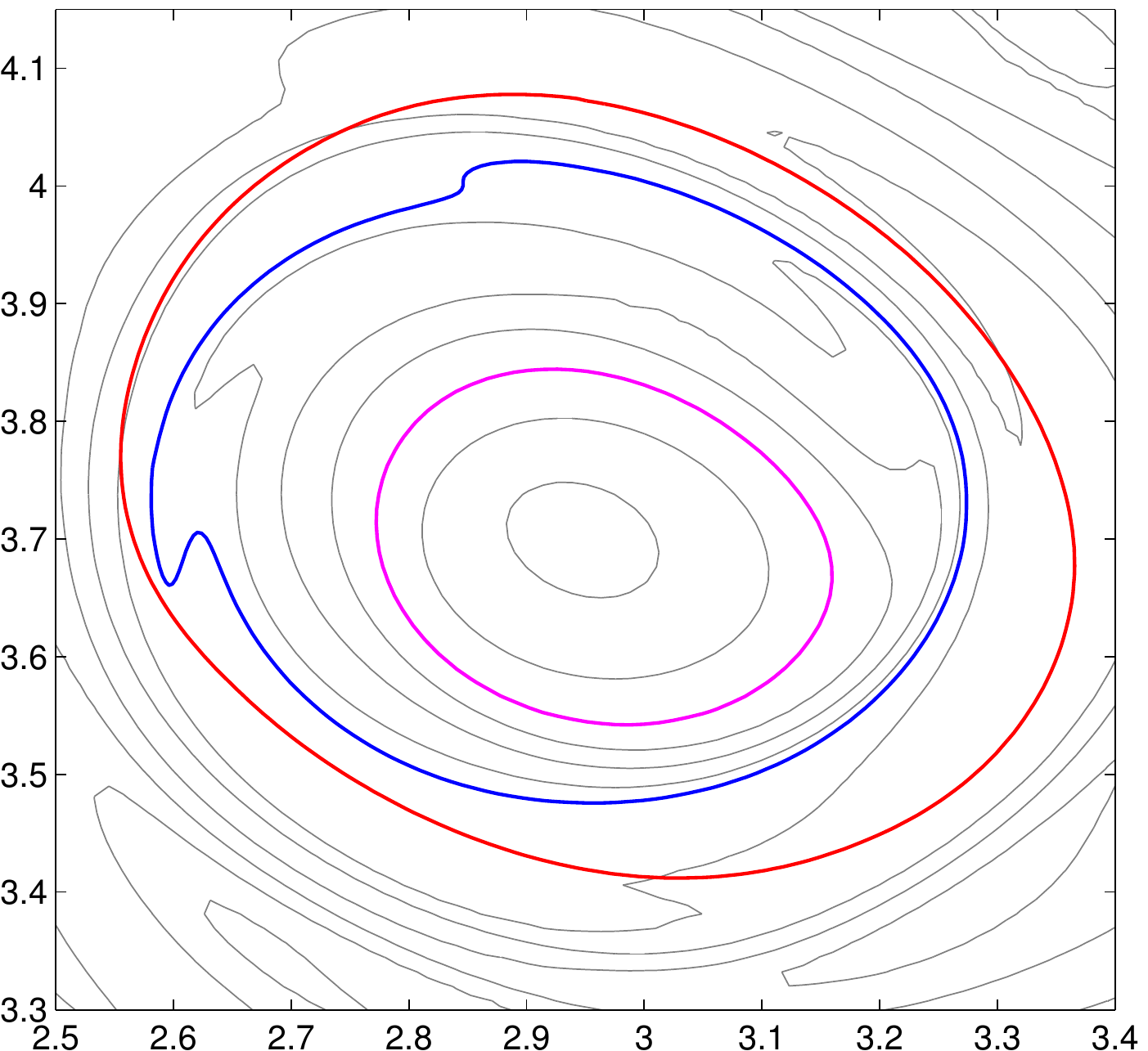}\hspace{.05\textwidth}
\includegraphics[width=0.45\textwidth]{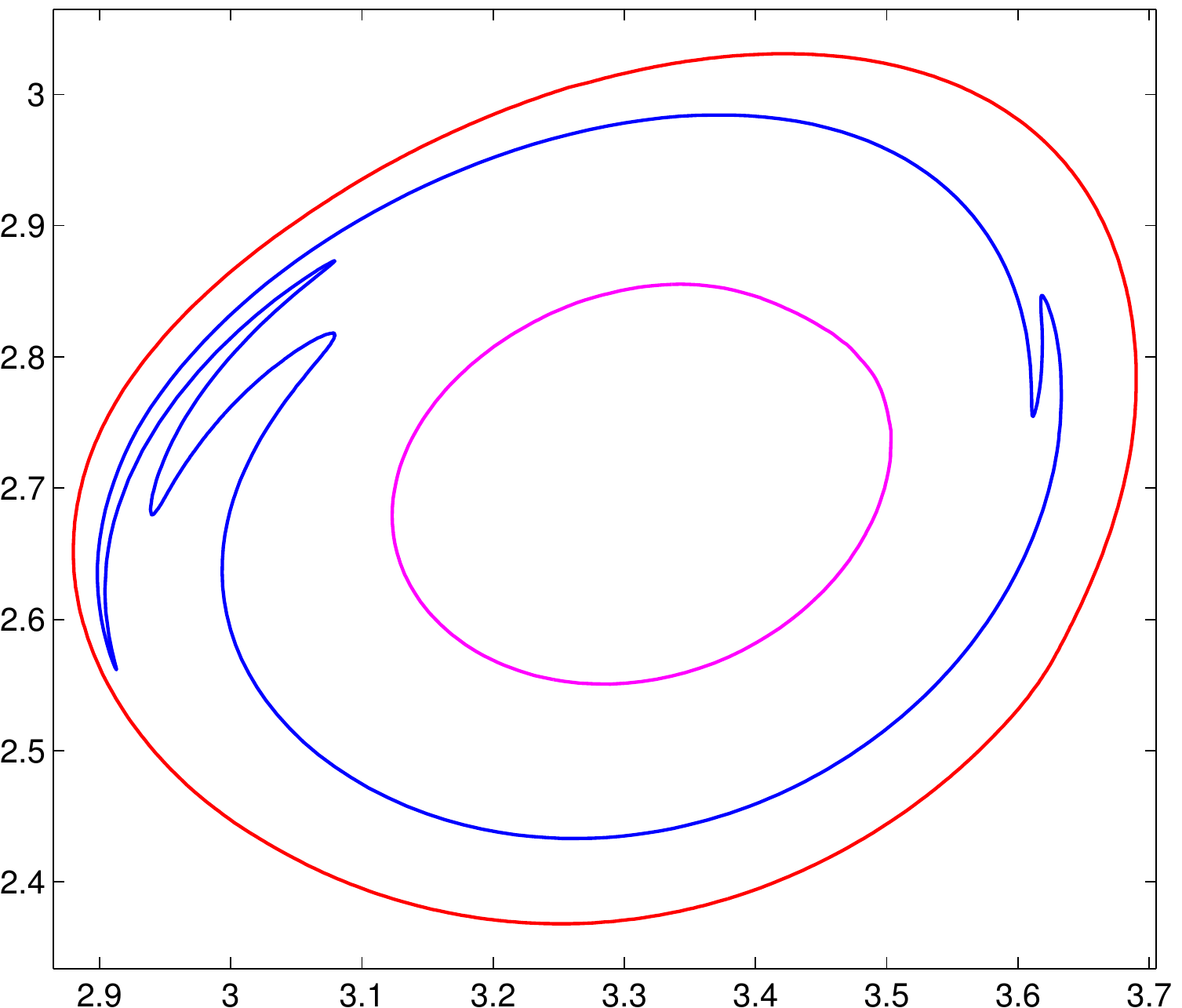}}
\caption{(a) Left: Vorticity contours (gray) and the Lagrangian vortex boundary
(red) for vortex 1 at time $t=0$. The blue curve marks the closed
vorticity contour that lays entirely inside the elliptic LCS. This
contour corresponds to $\omega=-0.3$. The magenta curve marks the
closest axisymmetric vorticity contour to the elliptic LCS. Right:
The Lagrangian vortex boundary and selected vorticity contours advected
to time $t=50$. (b) Same as (a) for vortex 3. The contour marked
by the blue curve corresponds to $\omega=-0.32$.}
\label{fig:cshr+vort}
\end{figure}

We find that the coherent vortex boundaries obtained
as outermost elliptic LCSs cannot be approximated by the instantaneous
vorticity contours at the initial time $t=0$. Figure \ref{fig:cshr+vort}
compares these vortex boundaries with the vorticity contours for two
of the coherent vortices. None of the vorticity contours approximates
the actual observed coherent vortex boundary of the Lagrangian frame.
In fact, the nearby vorticity contours are not axisymmetric, even
though that is intuitively expected for a coherent vortex boundary
\citep{mcwilliams_vort}. For instance, the closest vorticity
contour to the elliptic LCSs (blue curves in Fig. \ref{fig:cshr+vort}) notably lack axisymmetry. Their
advected positions at time $t=50$ develop filaments. In contrast,
the magenta-colored axisymmetric contours closest to the elliptic LCS
preserve their overall shape. These contours would, however, significantly
underestimate the true extent of the coherent fluid region.

{Similar observations can be made for the OW criterion.
The OW parameter 
\begin{equation}
Q=|S|^{2}-|\Omega|^{2},\label{eq:OW}
\end{equation}
measures instantaneous straining against instantaneous rotation. Here,
$S$ and $\Omega$ are, respectively, the symmetric and anti-symmetric
parts of the velocity gradient $\nabla u$. The matrix norms involved
are computed as $|S|^{2}=(\partial_{1}u_{1}-\partial_{2}u_{2})^{2}+(\partial_{1}u_{2}+\partial_{2}u_{1})^{2}$
and $|\Omega|^{2}=\omega^{2}$, where $(u_{1},u_{2})$ are the components
of the velocity field and $\omega=\partial_{1}u_{2}-\partial_{2}u_{1}$
is the vorticity field.}

{The subset of the domain where $Q>0$ is dominated
by strain, while $Q<0$ marks the regions dominated by vorticity.
Thus, the zero contour of this parameter encircling a vortex may be
expected to mark the outermost boundary of the vortical region. Several
authors have noted, however, that the zero contours of $Q$ will not
necessarily mark vortex-like structures (see, e.g., \citet{pasquero2001}
and \citet{koszalka2009}).}

{In practice, a negative-valued contour of $Q$ satisfying
$Q=-\alpha\sigma_{Q}$ is often considered as the vortex boundary\cite{pasquero2001},
where $\alpha$ is a positive constant and $\sigma_{Q}$ is the standard
deviation of the spatial distribution of $Q$. The constant $\alpha$
is somewhat arbitrary and must be tuned for a particular flow. \citet{pasquero2001,isern2006}
and \citet{Henson2008}, for instance, use $\alpha=0.2$ while \citet{koszalka2009}
use $\alpha=1$. }

\begin{figure}[t!]
{\centering}\subfloat[]{{\includegraphics[width=0.45\textwidth]{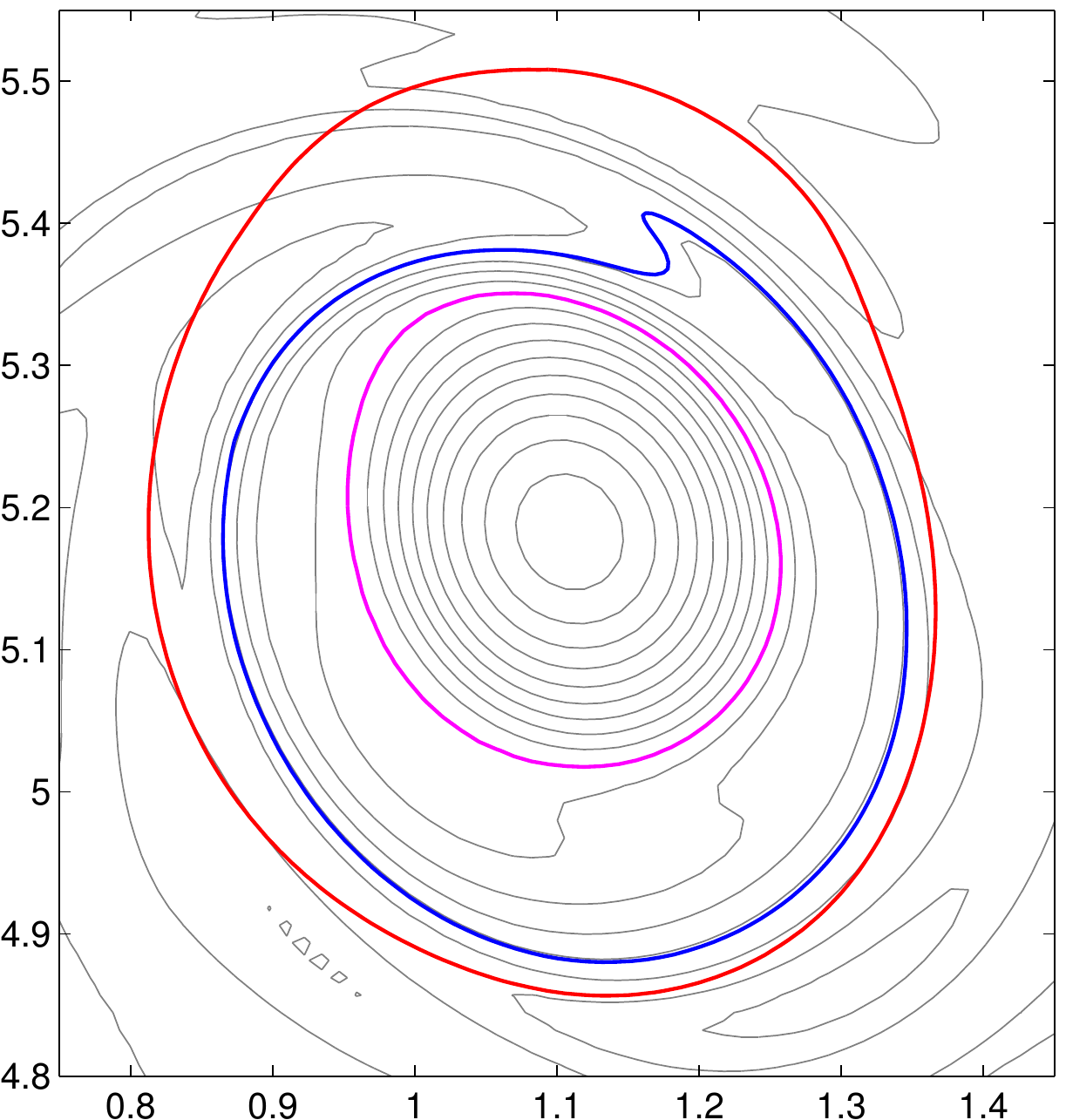}\hspace{.04\textwidth}
\includegraphics[width=0.5\textwidth]{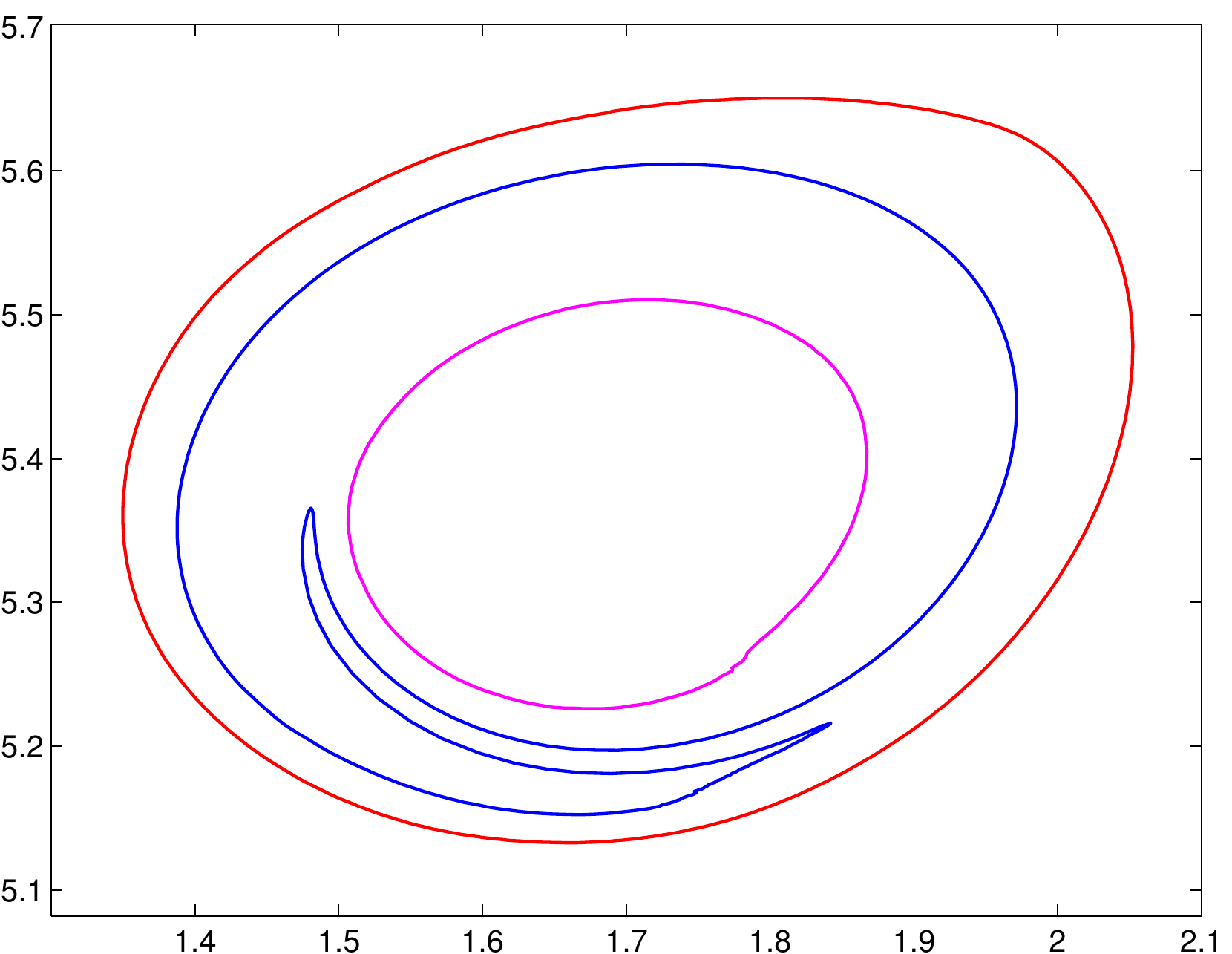}}}\\
\subfloat[]{{\includegraphics[width=0.4\textwidth]{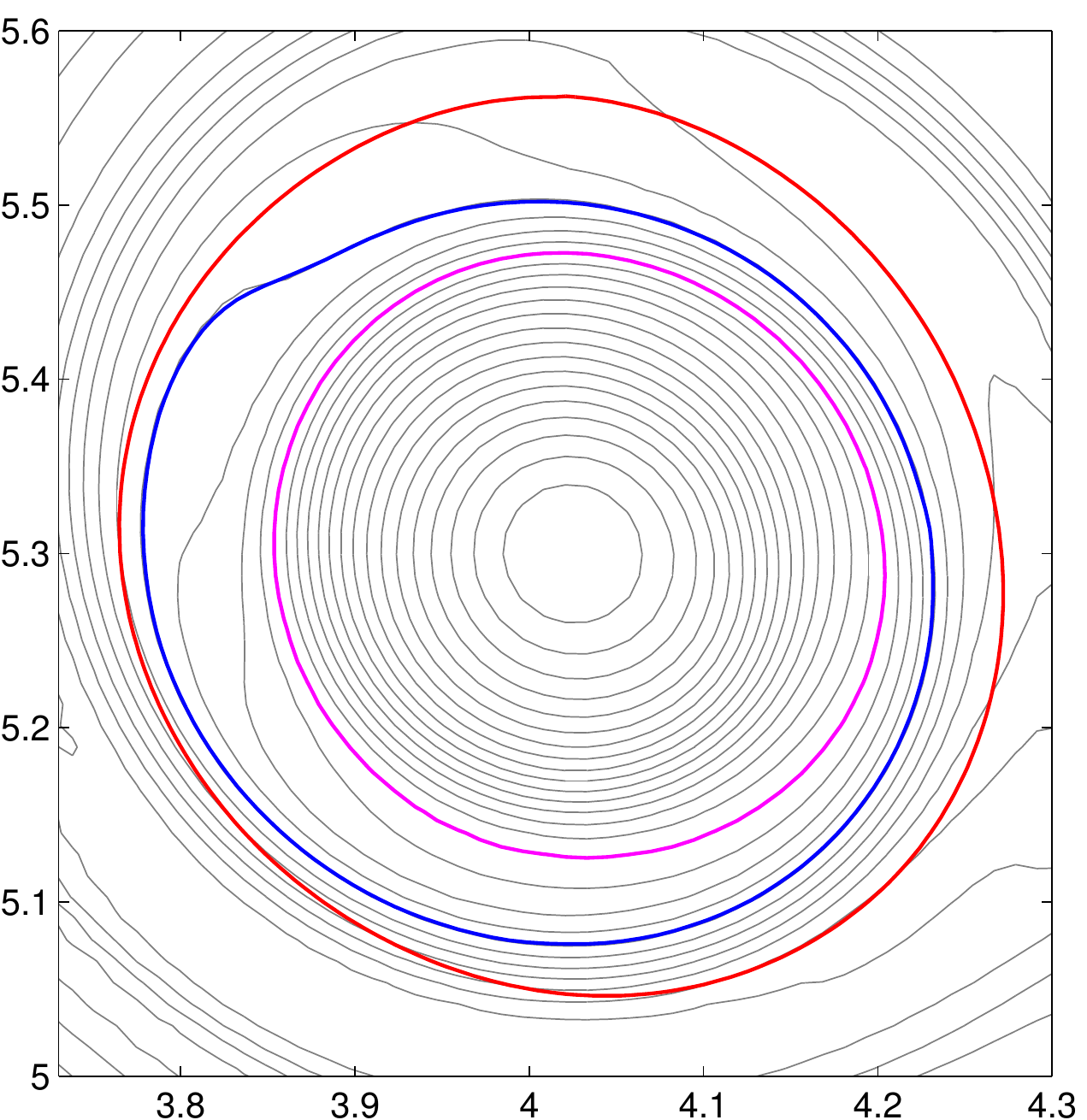}\hspace{.1\textwidth}
\includegraphics[width=0.41\textwidth]{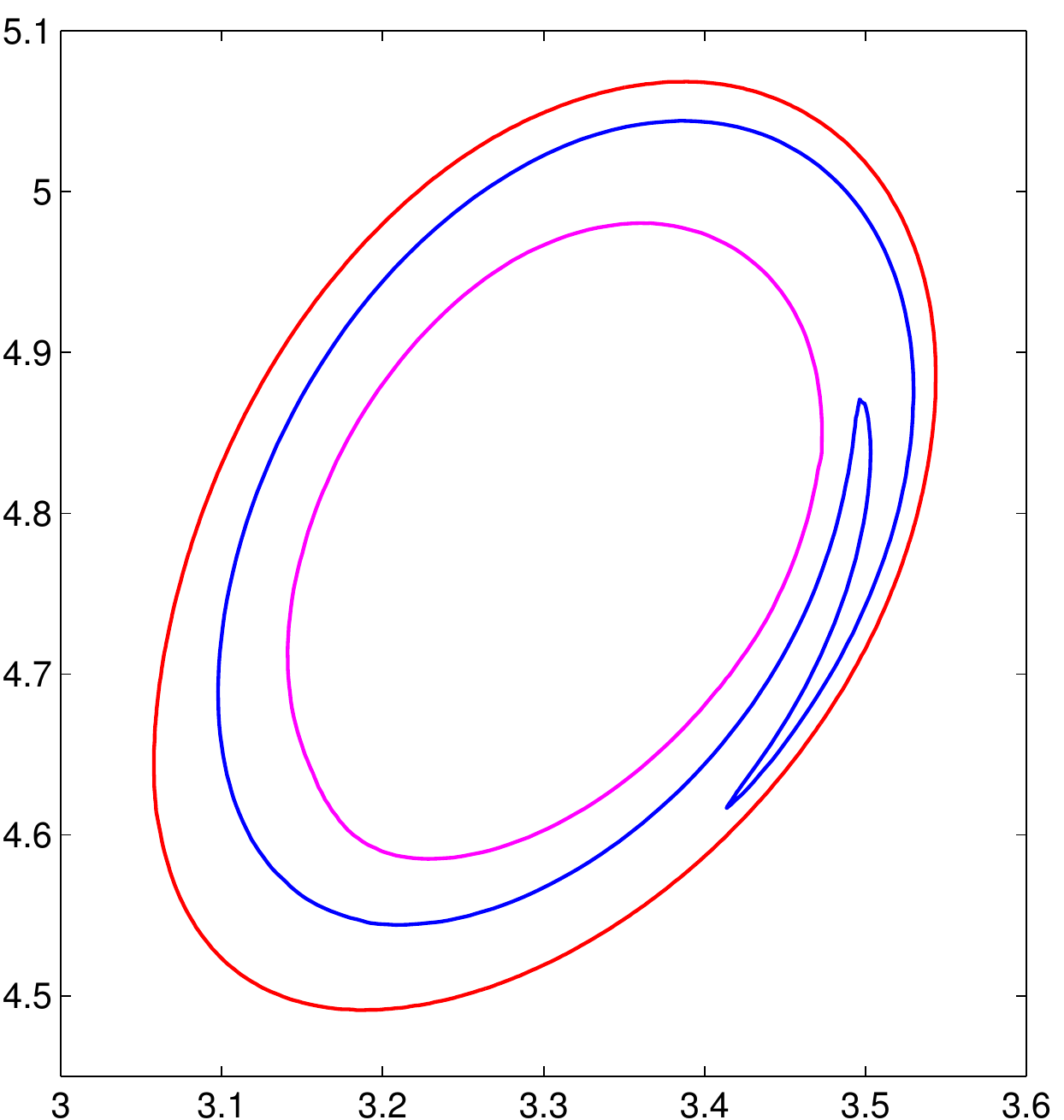}}}
\caption{(a) Left: OW contours (gray) and the Lagrangian vortex boundary (red)
for vortex 1 at time $t=0$. Two contours corresponding to $Q=-0.072\simeq-0.97\sigma_{Q}$
(blue) and $Q=-0.240\simeq-3.22\sigma_{Q}$ (magenta) are selected
for advection. Right: The Lagrangian vortex boundary and selected
OW contours advected to time $t=50$. (b) Same as (a) for vortex 2.
Here, the OW contours corresponding to $Q=-0.096\simeq-1.29\sigma_{Q}$
(blue) and $Q=-0.40\simeq-5.36\sigma_{Q}$ (magenta) are selected
for advection.}

{\label{fig:cshr+OW} }
\end{figure}

{A closer inspection of figure \ref{fig:cshr+OW}
reveals that none of the OW contours approximate well the true coherent
Lagrangian vortex boundary. The closest OW contour (blue curve) to
the outermost elliptic LCS lacks axisymmetry and develops substantial
filamentation under advection. The axisymmetric contour (magenta curve)
contained in the coherent vortex preserves its shape but seriously
underestimates the extent of the coherent region (as do axisymmetric
vorticity contours). This axisymmetric contour of the OW parameter
is also the outermost contour that remains in the $Q<0$ region over
the entire time interval $t\in[0,50]$.}

{We obtain similar conclusions about other OW-type
Eulerian indicators that have been developed to overcome the shortcomings
of the OW criterion (see, e.g., \citet{chong1990,tabor1994,kida1998,hua1998}).
\citet{hua1998}, for instance, consider the effect of higher-order
terms due to fluid acceleration. They arrive at the indicator parameters
$\lambda_{\pm}$ given by 
\[
\lambda_{\pm}=\frac{1}{4}Q\pm\frac{1}{2}\sqrt{|\dot{S}|^{2}-|\dot{\Omega}|^{2}},
\]
where $\dot{S}$ and $\dot{\Omega}$ denote, respectively, the instantaneous
rate of change of strain and vorticity along fluid trajectories. The
scalar $Q$ is the OW parameter, defined in \eqref{eq:OW}. The positive
extrema of $\lambda_{+}$ correspond to regions of instantaneously
strong stirring and dispersion. The negative extrema of $\lambda_{-}$,
on the other hand, mark the vortex regions. }
\begin{figure}[t!]
\centering \includegraphics[width=0.45\textwidth]{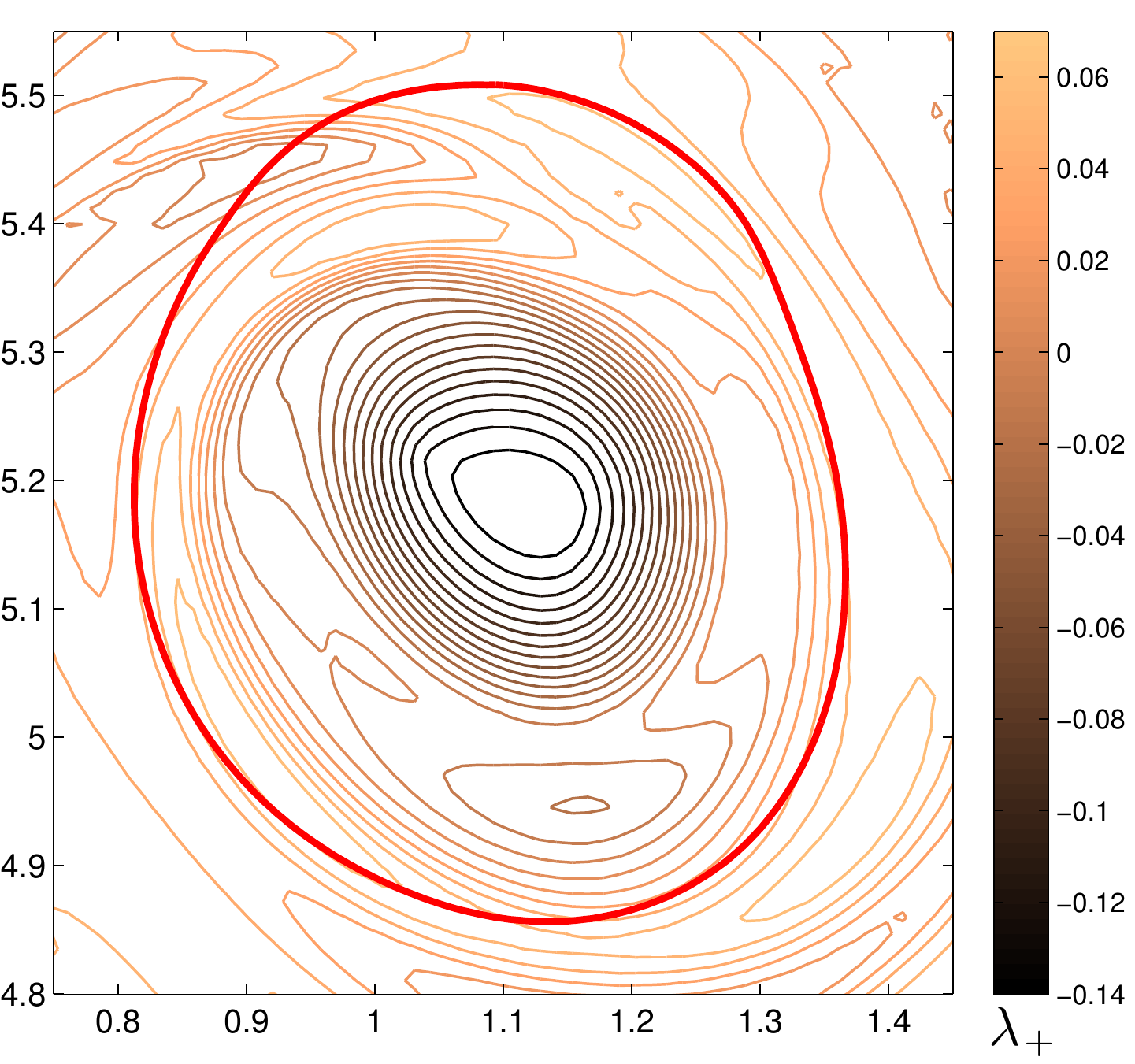}\hspace{.05\textwidth}
\includegraphics[width=0.45\textwidth]{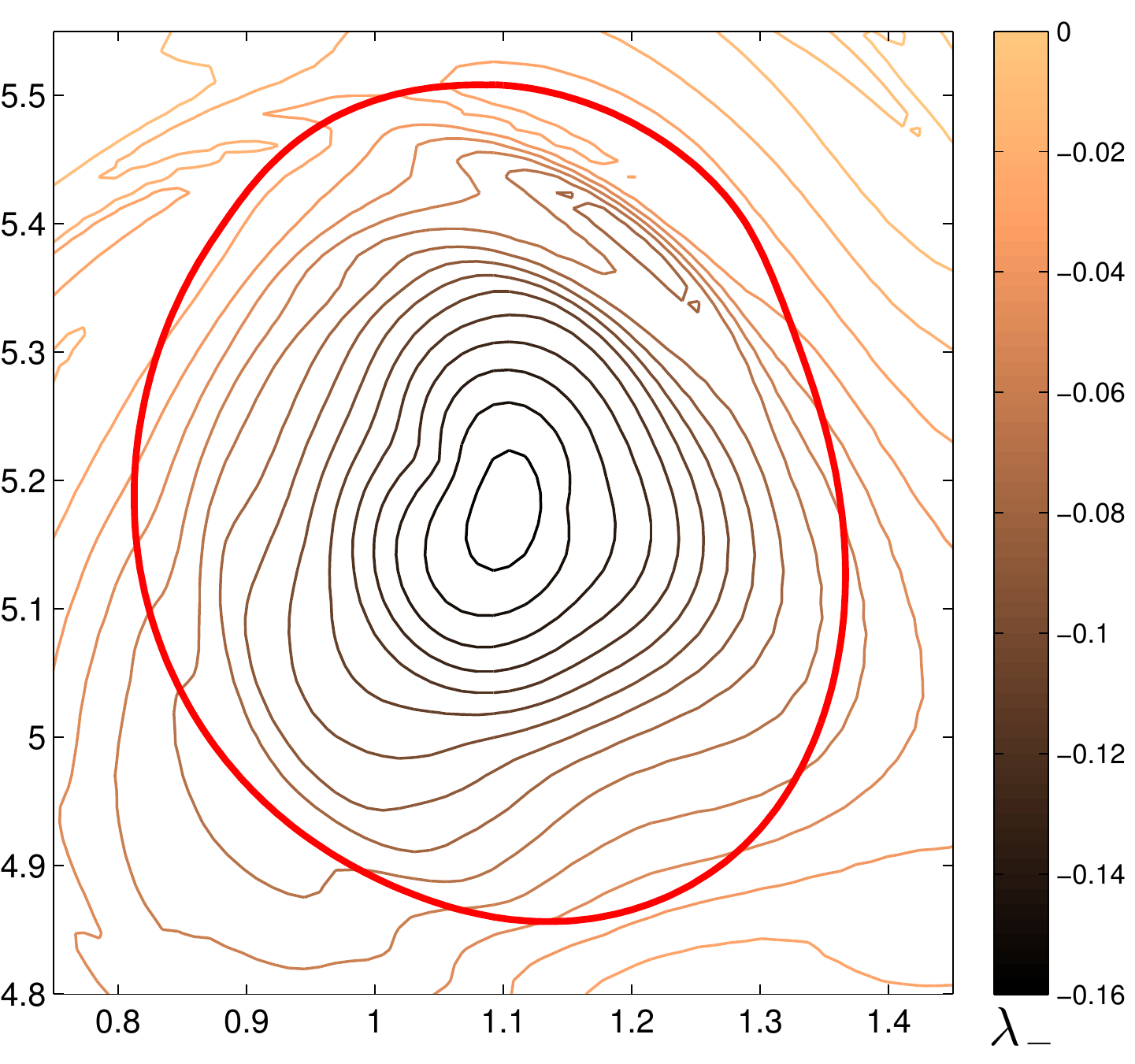}
\protect\caption{The contours of $\lambda_{+}$ (left) and $\lambda_{-}$ (right) around
vortex 1 at time $t=0$. The Lagrangian vortex boundary is shown with
thick red line.}
\label{fig:cshr+lambda}
\end{figure}

{As in the case of vorticity and the OW-parameter,
we find that the Lagrangian vortex boundaries cannot be inferred from
the contours of the $\lambda_{\pm}$ parameters (see figure \ref{fig:cshr+lambda}).
The axisymmetric contours of $\lambda_{\pm}$ remain coherent under
material advection over the time interval $t\in[0,50]$. They, however,
are significantly smaller (in enclosed surface area) than the true
Lagrangian vortex boundary marked by the elliptic LCS.}

\begin{figure}[t!]
\centering 
\includegraphics[width=0.6\textwidth]{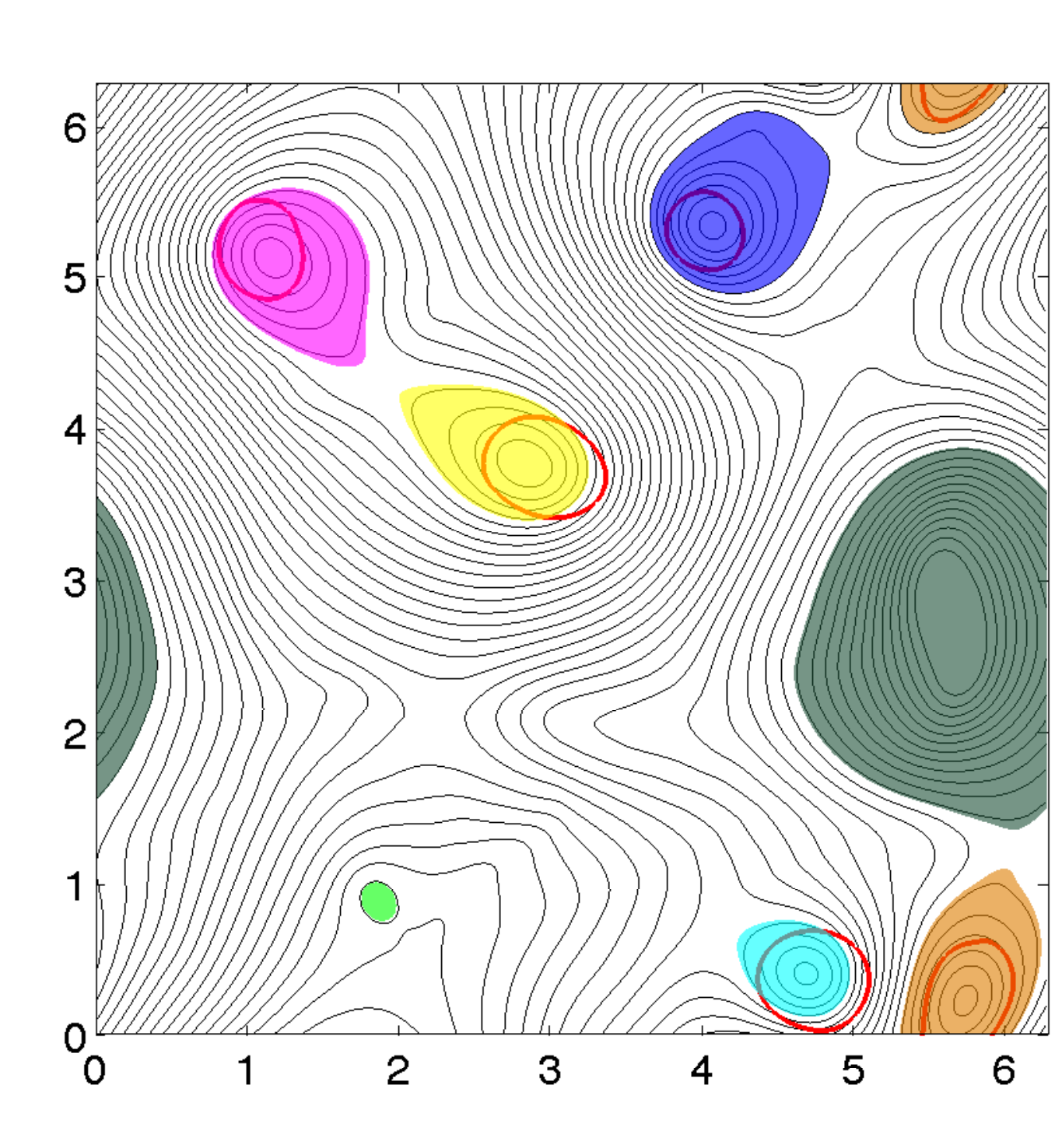}
\caption{Streamline-based eddies (colored patches), Lagrangian vortex boundaries
(red curves) and the streamlines (black curves) at the initial time
$t=0$.}
\label{fig:streamEddy}
\end{figure}

{The last Eulerian indicator we consider here is
the streamline-based eddy detection method proposed by \citet{entropy2Dturb}.
This method uses the topography of the instantaneous stream function
$\psi=-\Delta^{-1}\omega$ to locate a vortex. Specifically, a streamline-based
vortex boundary is locally the largest closed, numerically computed
streamline that does not enclose a saddle point. Figure \ref{fig:streamEddy}
shows the eddies detected in this fashion. The approach can be automated
using a cellular automata algorithm (see \citet{entropy2Dturb}, Appendix
A).}

{As noted by \citet{entropy2Dturb}, the streamline-based
eddy detection method seeks regions where strong vortical structures
may exist. Thus, each eddy island may, in principle, contain more
that one vortex, as is indeed the case for the dark green island of
figure \ref{fig:streamEddy}. Obtained from the instantaneous stream
function, however, the detected eddies are not guaranteed to preserve
their shape under advection. For instance, the (dark and light) green
patches quickly filament under advection.}

{Interestingly, the streamline-based eddy detection
also misses parts of the vortical structures: some of the Lagrangian
coherent vortices are not completely contained in streamline-based
eddy regions. This is the case for the vortices intersecting the yellow-
and cyan-colored patches in figure \ref{fig:streamEddy}. Other Lagrangian
coherent vortices happen to be fully contained in the blue-, magenta-
and brown-colored patches.}

{Compared to the number of Eulerian criteria, there
are far fewer Lagrangian diagnostics developed for quantifying coherent
vortices. These include the finite-time Lyapunov exponent \citep{ottino89,mixingFTLE},
mesoellipticity \cite{mezic_meso}, relative coherent pairs \citep{froyland,bolt_relCohSets},
shape coherence \citep{bolt_shapeCoh} and the ergodic partition of
time-averaged observables \citep{budivsic2012}. Here we evaluate
the performance of two of these in coherent Lagrangian vortex detection:
finite-time Lyapunov exponents and mesoellipticity.}

{The finite-time Lyapunov exponent (FTLE) measures
the maximal local stretching of material lines. For any point $x_{a}\in U$,
the FTLE corresponding to a time interval $[a,b]$ is defined as 
\begin{equation}
\Lambda(x_{a})=\frac{1}{2(b-a)}\log(\lambda_{2}(x_{a})),
\end{equation}
where $\lambda_{2}$ is the larger eigenvalue of the Cauchy--Green
strain tensor $C$. The FTLE measures the maximum separation of nearby
initial conditions over $[a,b]$. Therefore, its higher values suggest
regions of high stretching, and lower values generally indicate moderate
stretching. One envisions that low-FTLE regions to coincide with the
coherent Lagrangian vortex regions identified from our analysis.}

\begin{figure}[t!]
\centering\subfloat[]{\includegraphics[width=0.45\textwidth]{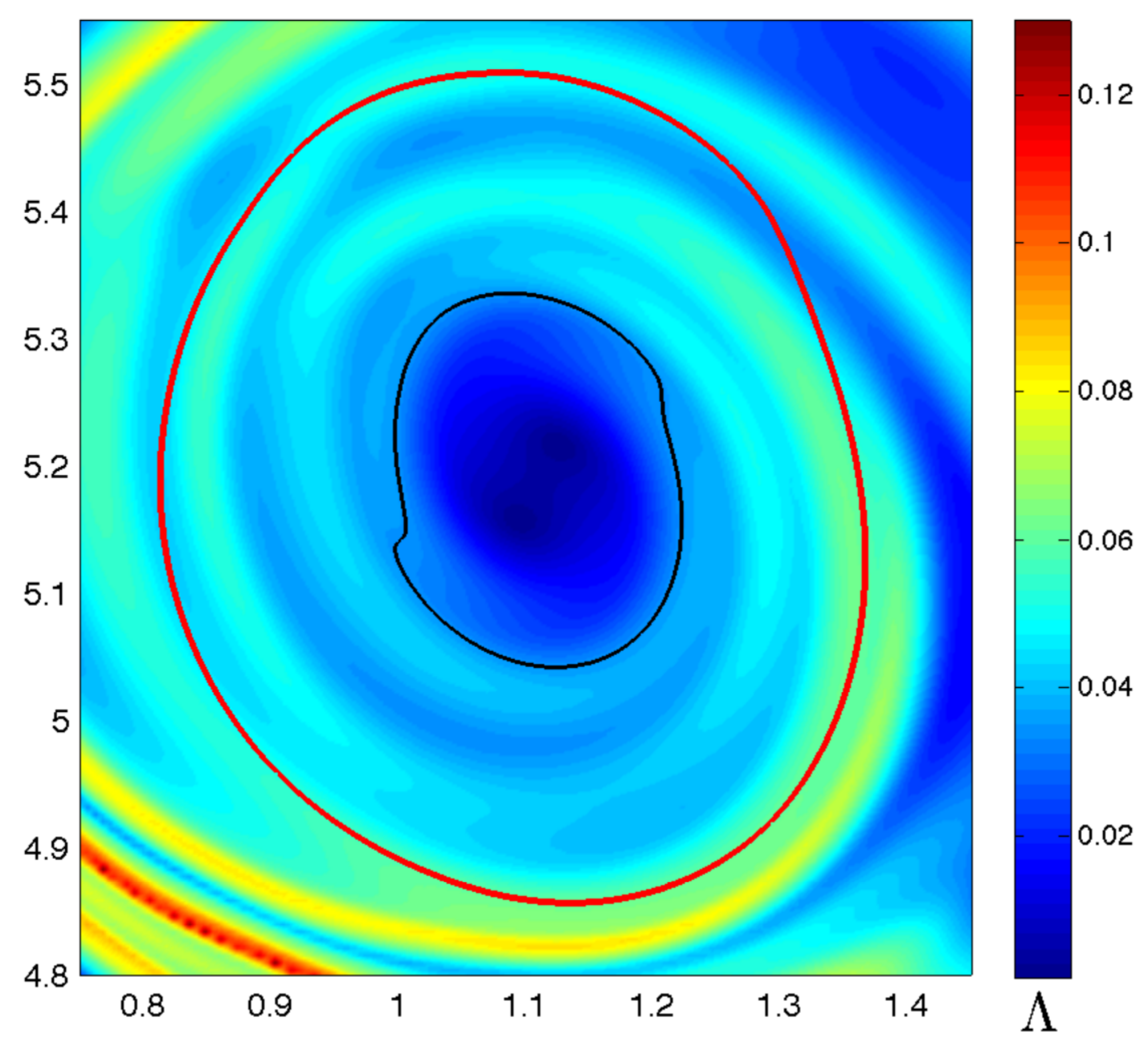}}
\hspace{.08\textwidth}
\subfloat[]{\includegraphics[width=0.43\textwidth]{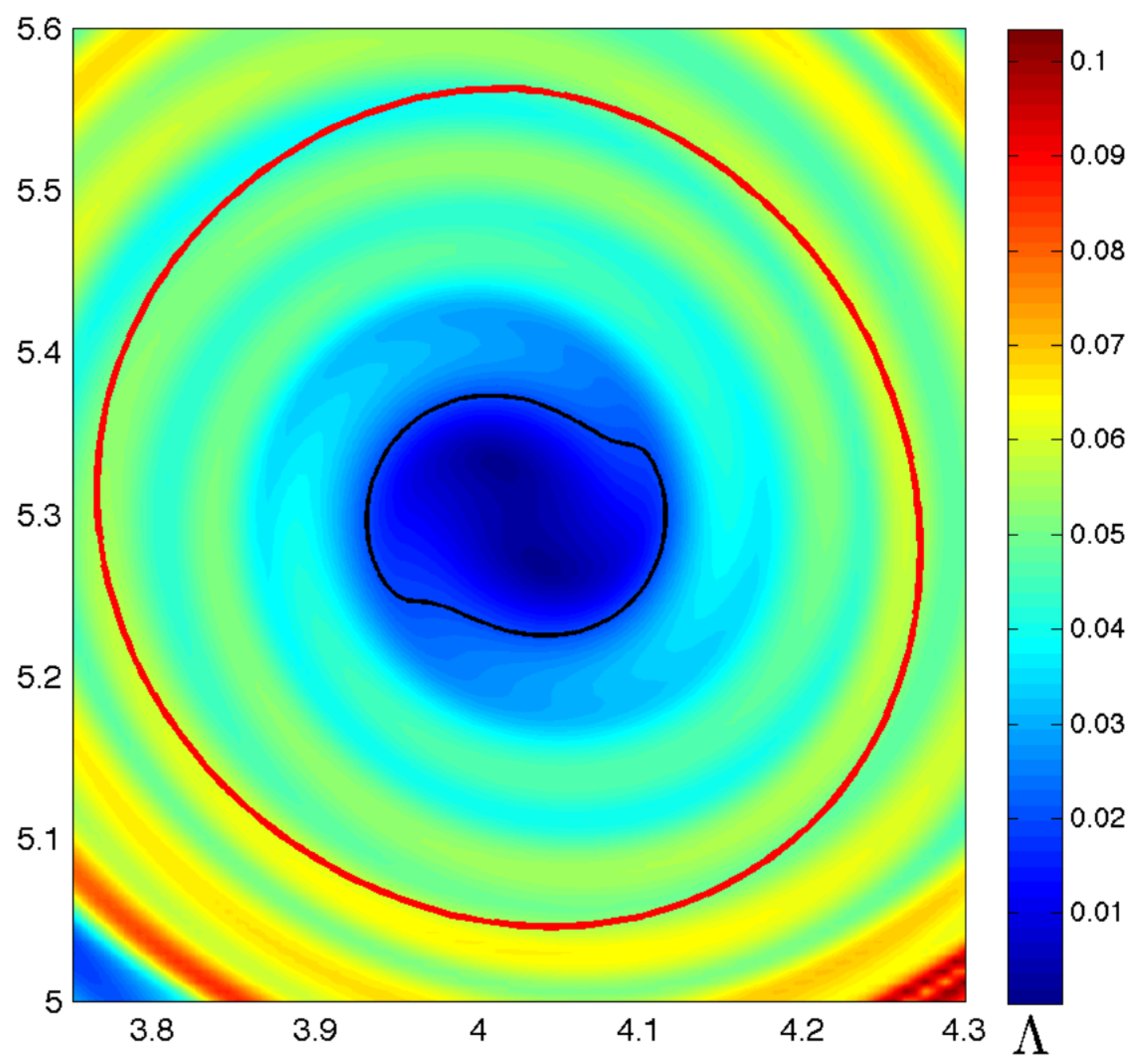}}
\caption{(a) Time $t=0$ position of the Lagrangian vortex boundary (red) for
vortex 1. The background color shows the FTLE field. The black curve
marks the FTLE contour with $\Lambda=3.45\times10^{-2}$. The FTLE
value is chosen such that the corresponding contour is the outermost,
almost-axisymmetric contour encircling the vortex core. (b) Same as
(a) for vortex 2. Here, the value of the FTLE contour is $\Lambda=2.0\times10^{-2}$}
\label{fig:cshr+ftle}
\end{figure}

{Figure \ref{fig:cshr+ftle} shows color-coded FTLE
values for vortices 1 and 2. Clearly, the Lagrangian vortex boundary
(red curves) cannot be inferred from the FTLE plot. In fact, locally
maximal values of FTLE spiral into the Lagrangian vortex boundary,
giving the wrong impression that it will stretch significantly under
advection.}

{In addition, FTLE contours around the vortex core
lack axisymmetry. The outermost, almost-axisymmetric FTLE contours
encircling the vortex cores (black curves) are still far from the
true vortex boundary marked by the elliptic LCS.}

{Now we consider a comparison between elliptic LCSs
and elliptic regions obtained from the Lagrangian mixing diagnostic
of \citet{mezic_meso}. This diagnostic classifies a trajectory starting
from a point $x_{a}$ as }{\emph{mesoellipitic}}{{}
in an incompressible flow, if the eigenvalues of the deformation gradient
$DF(x_{a})$ lie on the complex unit circle. Mesoelliptic trajectories
are expected to lie in a vortical region. In contrast, if the eigenvalues
of $DF(x_{a})$ are off the complex unit circle , the trajectory is
classified as }{\emph{mesohyperbolic}}{{}
and is expected to lie in a strain-dominated region.}

\begin{figure}[t!]
\includegraphics[width=0.8\textwidth]{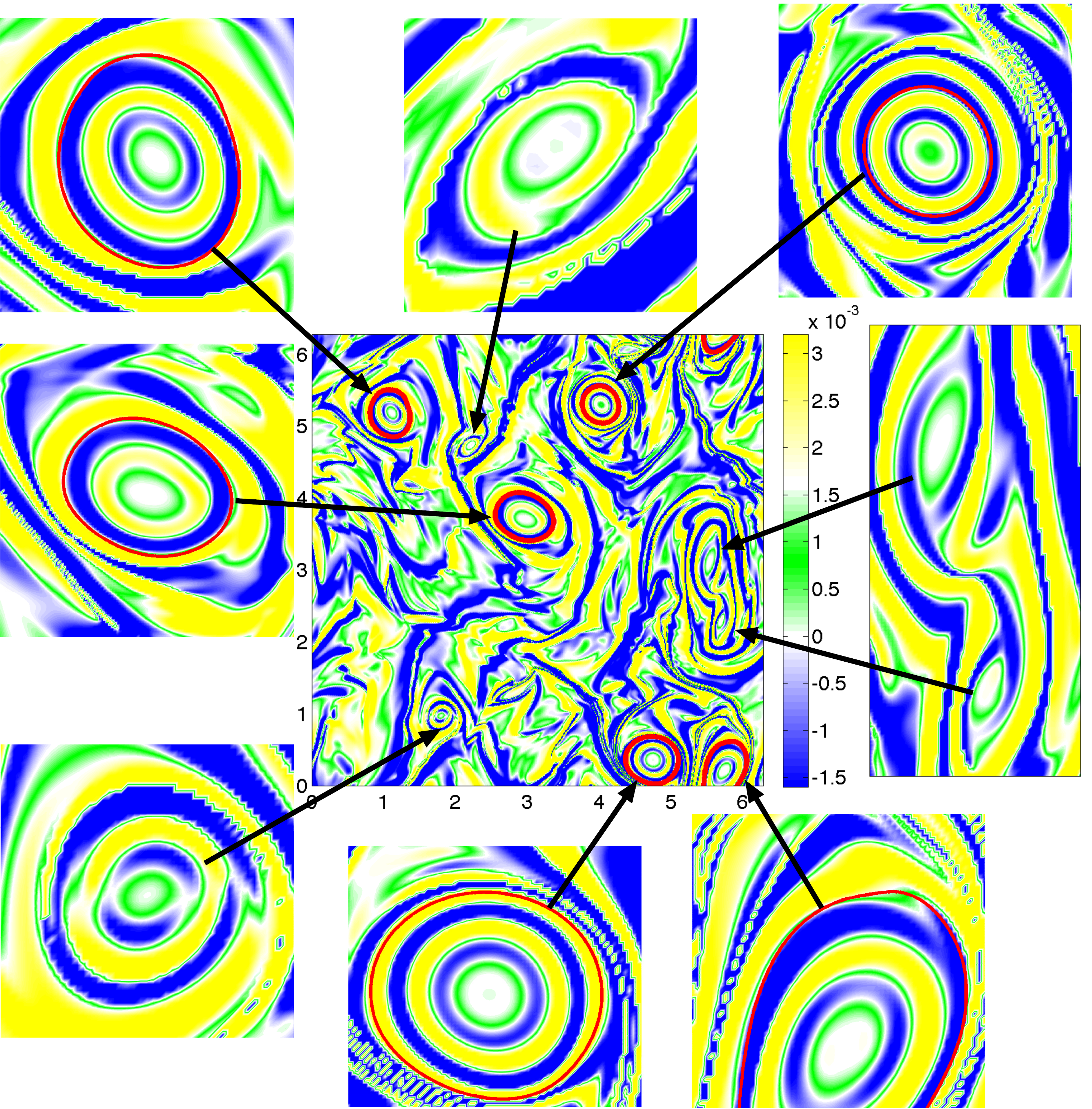}
\protect\caption{Hypergraph map for the turbulent flow marking mesohyperbolic (yellow
and blue color) and mesoelliptic (green and white) regions. The outermost
elliptic LCSs are shown as red curves. The three magnified regions,
where elliptic LCSs are absent, show examples of mesoelliptic regions
undergoing large stretching.}
\label{fig:meso}
\end{figure}

{Figure \ref{fig:meso} shows the hypergraph map
\cite{mezic_meso} for our turbulent flow, marking mesoelliptic (green
and white) and mesohyperbolic (yellow and blue) regions. As a rule,
the actual Lagrangian coherent vortex boundaries (i.e., the outermost
elliptic LCSs marked in red) always turn out to fall near the boundary
of a mesoelliptic (blue) annulus. Similar mesoellipotic annuli regions
exist, however, both inside and outside the Lagrangian vortex, thus
an a priori identification of the coherent vortex boundary cannot
be achieved based on mesoelliptic regions. Furthermore, substantial
portions of the actual Lagrangian vortices are diagnosed as mesohyperbolic
(annular yellow and blue regions). Likewise, a number of mesoelliptic
regions appear in non-coherent, hyperbolic mixing domains (compare
with figure \ref{fig:hypLCS}). These mesoelliptic areas undergo substantial
stretching and filamentation, and hence are false positives in coherent
vortex detection. We conclude that a systematic, a priori identification
of coherent Lagrangian vortex boundaries from mesoellipticity is not
possible in our turbulent flow.}

{}

\section{Conclusions}

{\label{sec:conclusion} We have used the variational
theory of \citet{bhEddy} to detect coherent material vortices in
a direct numerical simulation of two-dimensional Navier--Stokes turbulence.
We demonstrated that the vortex boundaries so obtained are optimal
in the sense that they are the outermost material lines enclosing
a vortex and retaining their shapes over long time intervals. They
are also frame-independent, threshold-free and Lagrangian by construction.}

{A comparison with other Eulerian methods (vorticity
contours, the Okubo--Weiss criterion, and the Hua--Klein criterion)
shows that the size of coherent material vortices in turbulence is
substantially larger than previously inferred from Eulerian indicators.
At the same time, the number of coherent vortices is lower than what
is signaled by the same indicators. This is consistent with the findings
in \citet{javier_agulhas}, who observed a similar trend for ocean
eddies in satellite-altimetry-based velocity fields of the South Atlantic.
We find that the superfluous vortices suggested by Eulerian indicators
are destroyed relatively quickly by the straining induced by repelling
and attracting Lagrangian coherent structures present in the flow.}

{We also compared our results with two Lagrangian
indicators: the finite-time Lyapunov exponent (FTLE) and the mesoellipticity
diagnostic of \citet{mezic_meso}. Low FTLE values generally indicate
the approximate position of vortex cores but do not provide an indication
of coherent Lagrangian vortex boundaries. Furthermore, FTLE lows also
occur in incoherent vortical regions as well, thus its use leads to
false positives in coherent mateal vortex detection. }

{As a rule, mesoelliptic annuli tend to form near
the coherent Lagrangian vortex boundaries. However, such annuli also
form both inside and outside coherent vortices, as well as in hyperbolic
mixing regions. Therefore, an
accurate and a priori identification of coherent Lagrangian vortices from mesoellipticity was not possible.}

{Compared to instantaneous Eulerian indicators, such
as Okubo-Weiss criterion, Lagrangian vortex detection is clearly computationally
more expensive. It requires accurate advection of a large ensemble
of fluid particles, as well as closed orbit detection in the vector
fields (\ref{eq:eta}). Therefore, developing cost effective computational
algorithms while staying faithful to the underlying theory is of great
interest (see \citet{leung2011,shadden2012,peikert2014}, for recent
developments).}

{Future theoretical work will focus on the correlation
between Lagrangian coherent vortices and the dynamical properties
of the flow, e.g., the scale-by-scale transfer of energy and enstrophy
\citep{lcs_Eflux}.}

{The streamline-based eddy detection discussed in
Section III.D was originally developed to study the formation of coherent
structures in decaying two-dimensional turbulence. Servidio et al.
\citet{entropy2Dturb} show that the local relaxation to coherent
structures is tied to the system's tendency to maximize a local notion
of entropy. It is of interest to re-examine this observation by adopting
the more accurate description of coherent vortices used in the present
paper.}
\begin{acknowledgments}
{We would like to thank Jan Feys for his help with
the implementation of the Fourier transform in MATLAB. We also thank
Daniel Blazevski for sharing his code for generating the hypergraph
map of figure \ref{fig:meso}.}
\end{acknowledgments}

\bibliographystyle{unsrtnat}
\bibliography{bibliog}

\end{document}